\documentclass[11pt]{article}
\usepackage[top=1in, bottom=1in, left=1in, right=1in]{geometry}
\usepackage{tikz,graphicx,amssymb,epstopdf,float,amsmath,amsthm,bm,color,epsfig,enumerate,enumitem,caption,hyperref,booktabs,multirow,array,comment,mathtools}
\usetikzlibrary{shapes,arrows,matrix,patterns,positioning}
\usepackage[norelsize,boxed,vlined,algo2e]{algorithm2e}
\usepackage{lineno}

\newtheorem{theorem}{Theorem}[section]

\newtheorem{conjecture}[theorem]{Conjecture}

\renewcommand{\appendix}[1]{
\section*{Appendix: #1}
}

\renewcommand{\O}[1]{\mathcal{O}\left(#1\right)}
\renewcommand{\L}[1]{\log\left(#1\right)}
\newcommand{\Ltwo}[1]{\log^2\left(#1\right)}
\newcommand{\Lthree}[1]{\log^3\left(#1\right)}

\newcommand{\bbC}{\mathbb{C}}
\newcommand{\bbR}{\mathbb{R}}

\renewcommand{\j}{j}

\newcommand{\rID}{{\it rID}}
\newcommand{\cID}{{\it cID}}

\newcommand{\sz}{\text{size}}
\newcommand{\svd}{\text{svd}}
\newcommand{\qd}{\text{;\quad}}
\SetKwProg{Fn}{Function}{ }{ }
\SetKwFunction{lrf}{LowRankFactorization}
\SetKwFunction{rsvd}{randomizedSVD}
\SetKwFunction{rp}{randperm}
\SetKwFunction{qr}{qr}
\SetKwFunction{lq}{lq}
\DeclarePairedDelimiter\ceil{\lceil}{\rceil}
\DeclarePairedDelimiter\floor{\lfloor}{\rfloor}

\makeatletter
\newcommand*{\extendadd}{
  \mathbin{
    \mathpalette\extend@add{}
  }
}
\newcommand*{\extend@add}[2]{
  \ooalign{
    $\m@th#1\leftrightarrow$%
    \vphantom{$\m@th#1\updownarrow$}
    \cr
    \hfil$\m@th#1\updownarrow$\hfil
  }
}
\makeatother

\begin{document}


\title{Rapid Application of the Spherical Harmonic Transform via Interpolative Decomposition Butterfly Factorization}

\author{James Bremer \\ Department of Mathematics, University of California, Davis, California, USA
    \and Ze Chen \\ Department of Mathematics, National University of Singapore, Singapore
    \and Haizhao Yang \\ Department of Mathematics, Purdue University, USA}

\maketitle

\begin{abstract}
We describe an algorithm for the application of the forward and inverse spherical harmonic transforms.  It is based on a 
new method for rapidly computing the forward and inverse associated Legendre transforms by hierarchically applying the 
interpolative decomposition butterfly factorization (IDBF).  Experimental evidence suggests that the complexity of
our method --- including all necessary precomputations --- is  $\O{N^2 \Lthree{N}}$  in terms of both flops and memory, 
where $N$ is the order of the transform.  This is nearly asymptotically optimal.  Moreover, unlike existing algorithms which 
are asymptotically optimal or nearly so, the constants in the running time and memory costs
of our algorithm are  small enough to make it competitive 
with state-of-the-art  $\O{N^3}$ methods at relatively small values of $N$ (e.g., $N=1024$). 
Numerical results are provided to demonstrate the effectiveness and numerical stability of the new framework.
\end{abstract}



\section{Introduction} \label{sec:intro}

This paper is concerned with the efficient application
of the forward and inverse spherical harmonic transforms (SHT).
These transformations play an important role in many scientific computing applications, including in the 
fields of numerical weather prediction and climate modeling \cite{fastSHT, weather2, weather1}, 
and are significant components in many numerical algorithms. The forward SHT of degree $N$ maps the   
coefficients in the expansion
\begin{equation}
    \label{eqn:sht_re}
    f(\theta,\phi) = \sum_{k=0}^{2N-1} \sum_{m=-k}^k \beta_{k,m} \overline{P}_k^{|m|} (\cos(\theta)) e^{im\phi},
\end{equation}
where $\overline{P}_k^{m}(x)$ denotes the $L^2$ normalized associated Legendre function of order $m$ and degree $k$,
to the values of the expansion at a grid of 
discretization nodes formed from the tensor product of a $2N$-point Gauss-Legendre quadrature in the variable $x = \cos(\theta)$
and a $(4N-1)$-point trapezoidal quadrature rule in the variable $\phi$.
The inverse SHT is, of course, the mapping which takes 
the values of the function $f(\theta,\phi)$ at the discretization nodes  
to the coefficients in the expansion.
More explicitly, the expansion $f$ is represented  either via the coefficients 
in (\ref{eqn:sht_re}) or by its values at the set of points
\begin{equation}
\left\{ \left(\theta_l, \phi_j\right)\ \ :
\ \ l=0,1,\ldots,2N-1, \ \ j=0,1,\ldots,4N-2
\right\},
\label{eqn:grid}
\end{equation}
where
\begin{equation}
    -1 < \cos(\theta_0) < \cos(\theta_1) < \dots < \cos(\theta_{2N-2}) < \cos(\theta_{2N-1}) < 1
\end{equation}
are the nodes of the $2N$-point Gauss-Legendre quadrature rule and 
$\phi_0, \phi_1, \dots, \phi_{4N-3}, \phi_{4N-2}$ are the equispaced nodes on $(0,2\pi)$ given
by the formula
\begin{equation}
    \phi_j = \frac{2\pi (j+\frac{1}{2})}{4N-1}, 
    \text{ for } j = 0, 1, \dots, 4N-3, 4N-2.
\end{equation}
The forward SHT maps the coefficients in (\ref{eqn:sht_re})
to the values of the expansion at the discretization nodes
\eqref{eqn:grid},  and the inverse SHT takes the values of the expansion
at the discretization nodes to the coefficients.

If we let
\begin{equation}
\label{eqn:sht_lt}
g(m,\theta)= \sum_{k=|m|}^{2N-1} \beta_{k,m} \overline{P}_k^{|m|} (\cos(\theta)),
\end{equation}
then  \eqref{eqn:sht_re} can be written as
\begin{equation}
\label{eqn:sht_fft}
 f(\theta,\phi) =    \sum_{m=-2N+1}^{2N-1} g(m,\theta) e^{i m \phi}.
\end{equation}
From \eqref{eqn:sht_fft}, it is clear that given
the values of $g(m,\theta)$ for each $m=-2N+1,\ldots,2N+1$ and each $\theta_0,\ldots,\theta_{2N-1}$,
the  values of $f(\theta,\phi)$
at the discretization nodes  \eqref{eqn:grid} can be computed in $\O{N^2 \L{N}}$ operations
by applying the  fast Fourier transform
$\O{N}$ times.  Similarly, the inverse of this operation, which takes
the values of $f(\theta,\phi)$ to those of  $g(m,\theta)$, can be calculated in $\O{N^2 \L{N}}$ 
operations using $\O{N}$ fast Fourier transforms.

We will refer to the mapping which, for a fixed $m$, takes the 
coefficients in the expansion \eqref{eqn:sht_lt} to the values of 
$g(m,\theta)$ at the $\O{N}$ discretization nodes in $\theta$ 
as the forward associated Legendre transform (ALT).  The inverse mapping,
which takes the values of $g(m,\theta)$ to 
the coefficients in the expansion,
will be referred to as the inverse ALT.
The naive approach to applying one of these transforms requires $\O{N^2}$ operations,
and using such an approach leads to an SHT with an  $\O{N^3}$ operation count.

There is a large literature devoted to accelerating the application of the associated Legendre transform
(we review it in Section~\ref{sec:relate}).  However, existing algorithms leave
much to be desired.  The most widely used class of methods allow for the application
of the ALT in $\O{N \log^\kappa\left(N\right)}$ operations, but only after
an $\O{N^2}$ precomputation phase.  Existing algorithms that have
quasilinear complexity (when all necessary precomputations are taken into account)
have such poor constants in their running times that they are slower than
the former class of methods at practical values of $N$.  Indeed, the current state-of-the-art 
method appears to be \cite{SHT3}, which has very favorable constants but requires a precomputation 
phase whose complexity is  $\O{N^2}$.

In this paper, we propose a new method for applying the forward
ALT whose total running time, including both the precomputation and application phases, 
is {$\O{(r(N))^2 N \Ltwo{N}}$},  
where $r(N)$ is a bound on the ranks of certain blocks
of the transformation matrix.  We conjecture that $(r(N))^2$ grows as
{$\O{\L{N}}$}.  Assuming this is correct, the total running time of our algorithm 
for applying the ALT
is  $\O{N \Lthree{N}}$. Proving a rigorous bound on $r(N)$ appears to be 
 quite difficult, and, for now, we are relying on experimental evidence regarding the running time of our algorithm.
Assuming our conjecture regarding the running time of our ALT is correct, the SHT
can be applied using our ALT  in  $\O{N^2 \Lthree{N}}$ time.

Our algorithm operates by  hierarchically applying the interpolative decomposition 
butterfly factorization (BF) \cite{IDBF,CHEN2020109427} (a newly proposed nearly linear 
scaling butterfly  algorithm \cite{BF, Butterfly1, Butterfly2,doi:10.1137/16M1074941,BF2}) and 
using  randomized low-rank approximation to speed up the calculation of the ALT.

\begin{figure}[!ht]
    \centering
    \begin{tabular}{ccc}
        \includegraphics[height=1.5in] {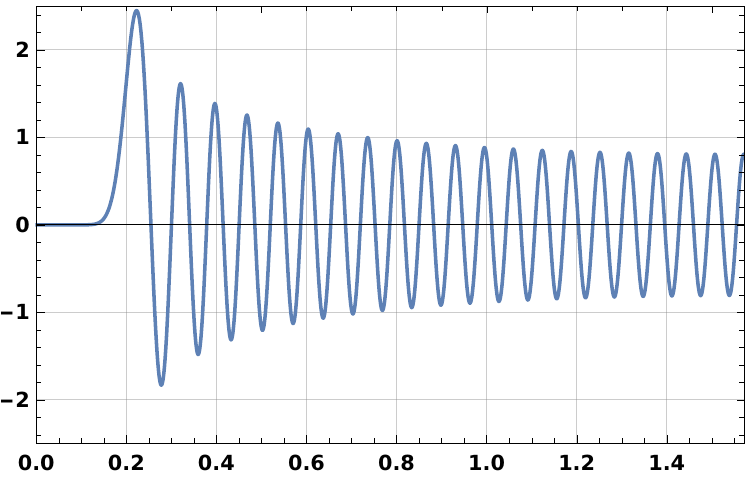} & & \includegraphics[height=1.5in] {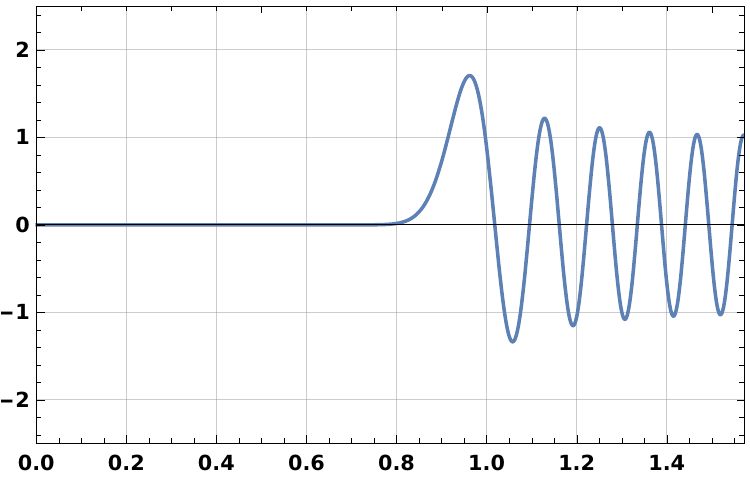}
    \end{tabular}
    \caption{On the left is a plot of the $L^2$ normalized associated Legendre function $\tilde{P}_n^m(\cos(\theta))$ 
in the case $n=100$ and $m=20$.  On the right is a plot of $\tilde{P}_n^m(\cos(\theta))$ when $n=100$ and $m=80$.}
    \label{fig:alf}
\end{figure}

Butterfly algorithms are a collection of techniques for rapidly applying the matrices
which result from discretizing oscillatory integral operators.
They exploit the fact that these matrices have the complementary low-rank property 
(see \cite{BF} for a definition).  A large class of special function transforms are integral operators of the appropriate type \cite{Butterfly2},
and consequently can be applied rapidly using butterfly algorithms.  Indeed, in the special case $m=0$,
the ALT  can be applied via standard butterfly algorithms in $\O{N \L{N}}$
time.  These results do not, however, extend to the case $m >0$.  In that event,
the associated Legendre functions are not oscillatory on the entire domain of interest.
Instead,  $\tilde{P}_n^m\left(\cos(\theta)\right)$ is nonoscillatory on the interval
\begin{equation}
\left(0,\arcsin{\left(\frac{\sqrt{m^2-1/4}}{n+1/2}\right)}\right)
\end{equation}
and oscillatory on 
\begin{equation}
\left(\arcsin{\left(\frac{\sqrt{m^2-1/4}}{n+1/2}\right)},\frac{\pi}{2}\right)
\end{equation}
(see Figure~\ref{fig:alf}, which contains  plots of $\tilde{P}_n^m(\cos(\theta))$ for two different pairs
of the parameters $n$ and $m$).
As a consequence, the integral operator associated with the 
ALT when $m >0$ is not of the purely oscillatory type whose discretizations
have the complementary low rank property.

\begin{figure}[!ht]
    \centering
    \begin{tabular}{ccc}
        \includegraphics[height=3in,trim={8cm 1cm 14.3cm 2cm},clip] {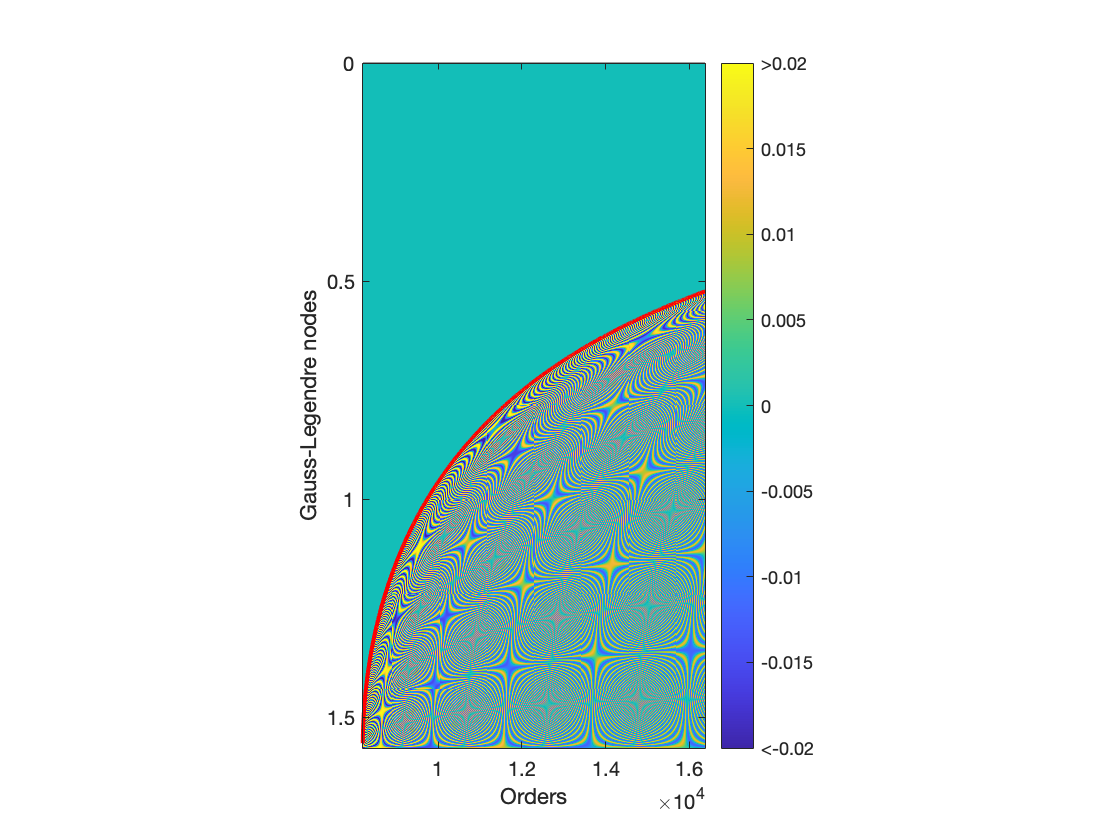} & & \includegraphics[height=3in,trim={8cm 1cm 14.3cm 2cm},clip] {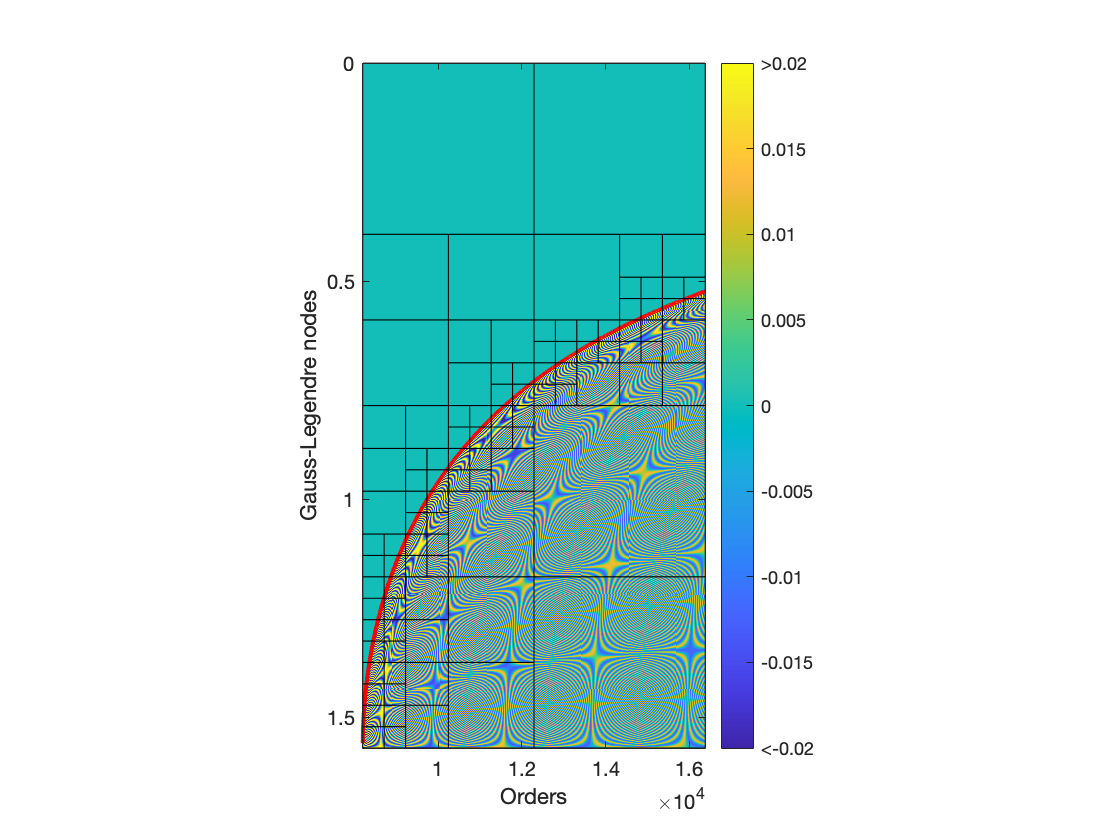} \\
        (a) & & (b)
    \end{tabular}
    \caption{An illustration of the partitioning process of an odd ALT matrix when $N = 8192$ and order $m = 8192$. 
(a) The odd matrix with a piecewise continuous curve (red color) indicating the positions of turning points. 
(b) The hierarchically partitioned blocks of the odd matrix.}
    \label{fig:lg}
\end{figure}

In order to overcome this difficulty we apply the following methodology:
\begin{itemize}
    \item We hierarchically partition the transformation matrix into purely oscillatory and purely non-oscillatory 
blocks (see Figure \ref{fig:lg} (b)).

    \item In the purely nonoscillatory blocks, the corresponding matrix is numerically low-rank, and hence its application 
to a vector can be accelerated to obtain linear scaling by randomized low-rank approximation algorithms.
    
    \item The matrices corresponding to purely oscillatory blocks admit complementary low-rank structures, 
the application of which to a vector can be accelerated via butterfly algorithms.  We use the relatively new interpolative decomposition butterfly factorization (IDBF) \cite{IDBF},
which yields nearly linear scaling in the degree $N$ of the ALT transform in both precomputation and application.

\end{itemize}

The scheme relies heavily on the algorithm of \cite{Bremer201815} for the numerical evaluation of the associated Legendre functions.  That algorithm allows each evaluation 
to be performed in time independent of the parameters $n$ and $m$.  If standard methods for calculating
$\tilde{P}_n^m$, which have running times
which grow with the parameters $n$ and $m$, were used instead, the running time of our algorithm 
for applying the ALT would no longer scale as $\O{N \Lthree{N}}$.

\subsection{Related works} \label{sec:relate}

There has been a significant amount of research aimed at accelerating the associated Legendre Transform
in order to more rapidly apply the spherical harmonic transform. 
In \cite{DRISCOLL1994202}, an algorithm for applying the ALT which results in an
SHT whose running time is $\O{N^2 \Ltwo{N}}$ is described.  However,
this algorithm suffers from increasing numerical instability as $N$ increases.  In \cite{fastSHT3} and 
\cite{SHT1}, asymptotically
optimal schemes for the ALT which are numerically stable are
described, but the constants in their running time make them unappealingly slow for practical values of  $N$.
The contribution \cite{SHT2} introduces a scheme based on the fast multiple method.
It can apply the SHT in $\O{N^2 \L{N}}$ operations after a
precomputation phase, the direct evaluation time of which is $\O{N^3}$ and could be reduced to nearly linear scaling in theory.  
However, this faster variant of the precomputation portion of the algorithm must be executed in extended precision arithmetic, 
which would most likely make it slow in applications.

The most widely-used algorithm today appears to be that of \cite{SHT3}.
It uses the butterfly transform described in  \cite{Butterfly2} to evaluate the ALT.
Each ALT takes $\O{N^2}$ and $\O{N \Lthree{N}}$ operations in the precomputation and application, respectively.  This, of course,
results in an SHT with a cost of $\O{N^3}$ for precomputation and $\O{N^2 \Lthree{N}}$ for application. 
A highly-optimized computational package based on \cite{SHT3} was developed in \cite{Seljebotn_2012}.  It is 
widely used
and most likely represents the current state-of-the-art for rapidly applying the SHT.
Though the application phase of this algorithm is nearly optimal, its precomputation phase is still
prohibitively expensive when $N$ is large. 

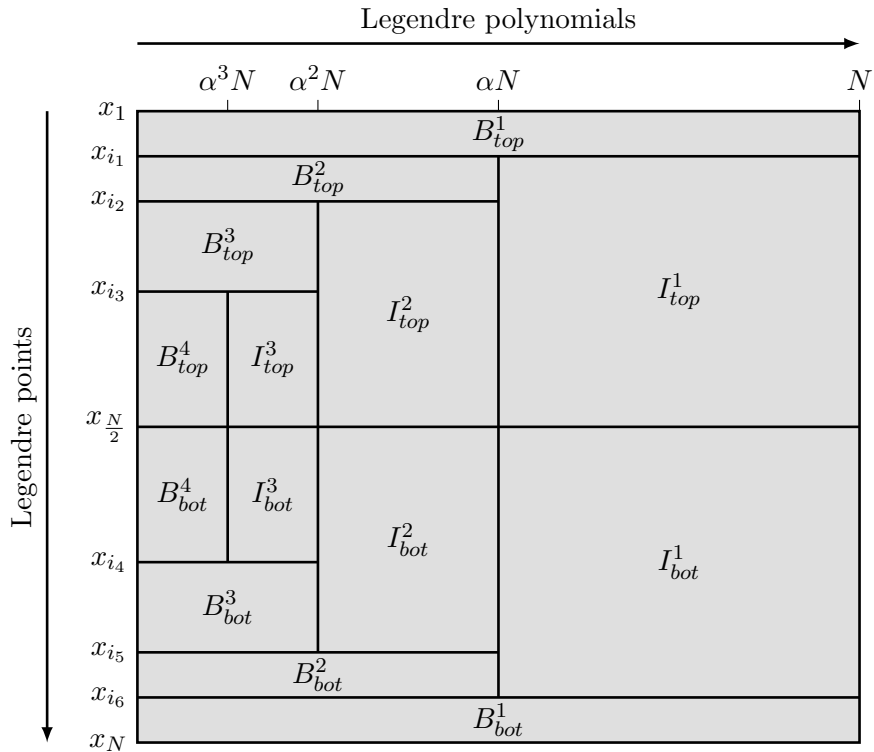
\begin{figure}[!ht]
\centering
\begin{tikzpicture}[scale=0.6]
\draw [line width=1pt, fill=gray!25] (0,0) -- (16,0) -- (16,14) -- (0,14) -- cycle;

\draw [line width=1pt] (0,1) -- (16,1);
\draw [line width=1pt] (0,7) -- (16,7);
\draw [line width=1pt] (0,13) -- (16,13);
\draw [line width=1pt] (0,2) -- (8,2);
\draw [line width=1pt] (0,12) -- (8,12);
\draw [line width=1pt] (0,4) -- (4,4);
\draw [line width=1pt] (0,10) -- (4,10);

\draw [line width=1pt] (8,1) -- (8,13);
\draw [line width=1pt] (4,2) -- (4,12);
\draw [line width=1pt] (2,4) -- (2,10);

\draw [-latex,line width=1pt] (-2,14) -> (-2,0);
\coordinate [label=left:\rotatebox{90}{Legendre points}] (X) at (-2,7);
\coordinate [label=left:$x_1$] (X) at (0,14);
\coordinate [label=left:$x_{i_1}$] (X) at (0,13);
\coordinate [label=left:$x_{i_2}$] (X) at (0,12);
\coordinate [label=left:$x_{i_3}$] (X) at (0,10);
\coordinate [label=left:$x_{\frac{N}{2}}$] (X) at (0,7);
\coordinate [label=left:$x_{i_4}$] (X) at (0,4);
\coordinate [label=left:$x_{i_5}$] (X) at (0,2);
\coordinate [label=left:$x_{i_6}$] (X) at (0,1);
\coordinate [label=left:$x_N$] (X) at (0,0);

\draw [-latex,line width=1pt] (0,15.5) -> (16,15.5);
\draw (2,14) -- (2,14.25);
\draw (4,14) -- (4,14.25);
\draw (8,14) -- (8,14.25);
\draw (16,14) -- (16,14.25);
\coordinate [label=above:Legendre polynomials] (X) at (8,15.5);
\coordinate [label=above:$\alpha^3 N$] (X) at (2,14.25);
\coordinate [label=above:$\alpha^2 N$] (X) at (4,14.25);
\coordinate [label=above:$\alpha N$] (X) at (8,14.25);
\coordinate [label=above:$N$] (X) at (16,14.25);

\coordinate (Bt1) at (8,13.5);
\node at (Bt1) {$B_{top}^1$};
\coordinate (Bb1) at (8,0.5);
\node at (Bb1) {$B_{bot}^1$};
\coordinate (Bt2) at (4,12.5);
\node at (Bt2) {$B_{top}^2$};
\coordinate (Bb2) at (4,1.5);
\node at (Bb2) {$B_{bot}^2$};
\coordinate (Bt3) at (2,11);
\node at (Bt3) {$B_{top}^3$};
\coordinate (Bb3) at (2,3);
\node at (Bb3) {$B_{bot}^3$};
\coordinate (Bt4) at (1,8.5);
\node at (Bt4) {$B_{top}^4$};
\coordinate (Bb4) at (1,5.5);
\node at (Bb4) {$B_{bot}^4$};

\coordinate (It1) at (12,10);
\node at (It1) {$I_{top}^1$};
\coordinate (Ib1) at (12,4);
\node at (Ib1) {$I_{bot}^1$};
\coordinate (It2) at (6,9.5);
\node at (It2) {$I_{top}^2$};
\coordinate (Ib2) at (6,4.5);
\node at (Ib2) {$I_{bot}^2$};
\coordinate (It3) at (3,8.5);
\node at (It3) {$I_{top}^3$};
\coordinate (Ib3) at (3,5.5);
\node at (Ib3) {$I_{bot}^3$};

\end{tikzpicture}

\caption{Block partitioning of the Legendre-Vandermonde matrix in \cite{fastSHT}, when $N = 1024$. $x_i$, $i = 1, 2, \dots, N$, are the Legendre points. The parameters $i_1, i_2, \dots, i_6$ and $\alpha$ are the computed partitioning coefficients, which are able to divide the Legendre-Vandermonde matrix into boundary parts (denoted by symbol $B$) and internal (denoted by symbol $I$) parts. The internal parts can be compressed by the BF while the boundary parts are directly computed for the corresponding matvecs.}
\label{fig:lt}

\end{figure}

In \cite{fastSHT} an algorithm for applying the ALT which bears some similarities to our scheme was proposed.
It  operates by  partitioning the transformation matrix in the  manner shown in 
Figure \ref{fig:lt}.     The application phase of the resulting algorithm has lower complexity
than that used in \cite{SHT3}
and yields somewhat improved accuracy (roughly an extra digit of precision).
However, the method of \cite{fastSHT} still requires a precomputation phase whose
running time behaves as $O(N^3)$.

In \cite{Slevinsky2019}, an algorithm that makes use of a rapid transformation between spherical harmonic expansions and bivariate Fourier series via the butterfly transform and hierarchically off-diagonal low-rank matrix decompositions. 
Although the application time of this algorithm $\O{N^2 \Ltwo{N}}$, 
it requires a precomputation whose running time grows as $\O{N^3 \L{N}}$.

\subsection{Outline of this paper}

The rest of the paper is organized as follows. 
In Section~\ref{sec:pre}, we discuss existing low-rank matrix factorization techniques. 
Section~\ref{sec:algo} proposes a new algorithm for applying the Legendre Transform
which is based on these factorization techniques.
In Section~\ref{sec:analysis}, we discuss the computational complexity of our algorithm.
Again, we do not have a rigorous bound on its running time, but we estimate
it under an assumption on the behavior of the ranks of certain subblocks of the 
matrix discretizing the ALT.
Section~\ref{sec:results} describes several numerical experiments conducted to assess
the efficiency of the proposed algorithm.

For simplicity, we adopt MATLAB notations for the algorithm described in this paper: given row and column 
index sets $I$ and $J$, $K(I,J)$ is the submatrix with entries from rows in $I$ and columns in $J$; the index set for 
an entire row or column is denoted as $``:"$.

\section{Low rank factorizations and butterfly algorithms}
\label{sec:pre}

In this section, we discuss several existing algorithms which exploit rank deficiency to rapidly apply certain
classes of matrices.  These are used by the algorithm of Section~\ref{sec:algo}
for the application of the ALT.
Subsection~\ref{sec:ID} outlines the linear scaling interpolative decomposition (ID) method, 
which is an important tool for the interpolative decomposition butterfly factorization (IDBF) discussed
in Subsection~\ref{sec:IDBF}. Subsection~\ref{sec:rSVD} describes low-rank approximation via randomized sampling.

\subsection{Linear scaling interpolation decomposition} \label{sec:ID}

This subsection reviews the algorithm of  \cite{IDBF} for the construction of interpolative
decompositions.

A column interpolative decomposition, which will abbreviate by $\cID$, of  $A\in \bbC^{m \times n}$  is a 
factorization of the form
\begin{equation}
 A \approx A(:,q)V,
\label{eqn:id}
\end{equation} 
where $q$ is an index set specifying $k$ columns of $A$ and $V$ is a $k \times n$ matrix.
The set  $q$ is called the {\it skeleton} index set, 
and the rest of the indices are called {\it redundant} indices.   The matrix $V$ is called the column
interpolation matrix.    
The algorithm described in this section takes as input a desired precision $\epsilon$ and adaptively 
determines $k$ such that
\begin{equation}
\left\| A - A(:,q)V\right\|_2 \leq \epsilon.
\end{equation} 
The numerical rank of $A$ to precision $\epsilon$ is defined via
\begin{equation}
    k_{\epsilon}=\min \left\{\operatorname{rank}(X): X \in \mathbb{C}^{m \times n},\|A-X\|_{2} \leq \epsilon\right\},
\end{equation}
and it is the optimal possible value of $k$.  In most cases, the algorithm of this section
forms factorizations with $k$ equal to or only slightly larger than $k_\epsilon$.

The algorithm also takes as  input
a parameter $r_k$, which we refer to as the ``adaptive rank,''
 which serves as an upper bound for
the rank of $A$. It proceeds by first constructing an index set $s$
containing $t \cdot r_k$ rows of $A$ chosen from the 
Mock-Chebyshev grids as in \cite{NUFFTorBF,Mock1,Mock2} or 
randomly sampled points.  Here, $t$ is an oversampling parameter.

We next compute a rank revealing QR decomposition of  $A(s,:)$.  That is,
we decompose $A(s,:)$ as 
\begin{equation}
\label{eq:pivotedQR}
    A(s,:)\Lambda \approx QR = Q[R_{1} \ R_{2}],
\end{equation}
where the columns of $Q\in \bbC^{tk \times {k}}$ are an orthonormal set in $\bbC^m$, $R\in \bbC^{{k}\times n}$ 
is upper trapezoidal, and $\Lambda \in \bbC^{n\times n}$ is a carefully chosen permutation matrix such that 
$R_{1}\in \bbC^{{k\times k}}$ is nonsingular.  The value of $k$ is chosen so that
the $L^2$ error in the approximation (\ref{eq:pivotedQR}) is somewhat smaller than $\epsilon$.
We now define
\begin{equation}
    T = (R_{1}(1:{k},1:{k}))^{-1}[R_1(1:{k},{k+1:tk}) \ R_{2}(1:{k},:)]\in \mathbb{C}^{{k\times (n-k)}},
\end{equation}
such that \[A(s,q)=QR_1(1:{k},1:{k}).\]   Then 
\begin{equation}
    A(s,:) \approx  A(s,q)V
\end{equation} 
with $V=[I ,T]\Lambda$ and the approximation error determined by the error in the rank-revealing QR decomposition.
Moreover, 
\begin{equation}
    A \approx A(:,q)V
\label{eqn:id2}
\end{equation} 
with an approximation error coming from the QR truncation and the error incurred in \eqref{eqn:id} when performing interpolation from the subsampled rows of $A$ using the interpolation matrix $V$.
When the obtained accuracy is insufficient, the procedure is repeated with an increased $k$.
Using the steps outlined above, the construction of this factorization requires $\O{nk^2}$ operations
and $\O{nk}$ storage.

A row interpolative  decomposition (abbreviated $\rID$) of the form

\begin{equation}
    A \approx U A(q,:)
\end{equation}
can be constricted in a similar fashion in $\O{mk^2}$ operations using $\O{mk}$ storage.  We refer to $U$
as the {\it row interpolation matrix}.

\subsection{Interpolative decomposition butterfly factorization}
\label{sec:IDBF}

In this section, we briefly discuss the properties of the interpolative decomposition butterfly factorization,
and the algorithm of \cite{IDBF} for producing it.   
We refer the reader to  \cite{IDBF} for a detailed discussion.

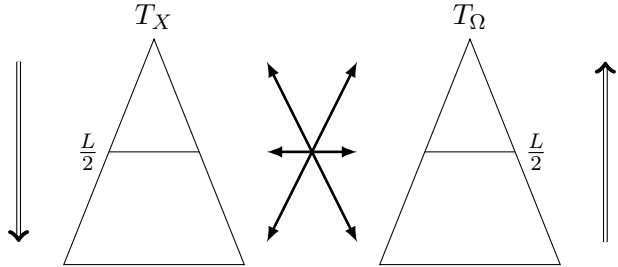
\begin{figure}[!ht]
\centering
\begin{tikzpicture}[scale=0.3]

\draw [->,double,double distance=1pt] (-6,9) -> (-6,1);

\draw (0,10) -- (-4,0) -- (4,0) -- cycle;

\draw [->,double,double distance=1pt] (20,1) -> (20,9);

\draw (14,10) -- (10,0) -- (18,0) -- cycle;

\draw [latex-latex,line width=1pt] (9,5) -> (5,5);
\draw [latex-latex,line width=1pt] (9,9) -> (5,1);
\draw [latex-latex,line width=1pt] (9,1) -> (5,9);

\coordinate [label=above:$T_X$] (X) at (0,10);
\coordinate [label=above:$T_\Omega$] (X) at (14,10);

\draw (-2,5) -- (2,5);
\draw (12,5) -- (16,5);

\coordinate [label=right:$\frac{L}{2}$] (X) at (16,5);
\coordinate [label=left:$\frac{L}{2}$] (X) at (-2,5);
\end{tikzpicture}
\caption{Trees of the row and column indices.
    Left: $T_X$ for the row indices $X$.
    Right: $T_\Omega$ for the column indices $\Omega$.
    The interaction between $A\in T_X$ and $B\in T_\Omega$
    starts at the root of $T_X$ and the leaves of $T_\Omega$. }
\label{fig:domain-tree-BF}
\end{figure}

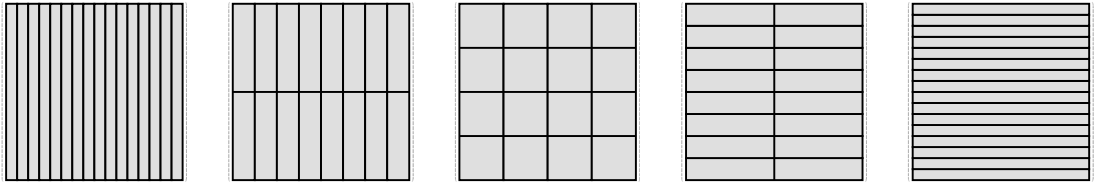
\begin{figure}[!ht]
\begin{minipage}{\textwidth}
\centering
\resizebox{2.5cm}{!}{
\begin{tikzpicture}[baseline=-0.5ex]
      \tikzset{every left delimiter/.style={xshift=-1ex},every right delimiter/.style={xshift=1ex}}
      \matrix (mat) [left delimiter=(, right delimiter=)] {
      \draw;
\foreach \j in {0,2,...,30}{
    \draw[fill=gray!25,line width=10pt] (\j,32) rectangle (\j+2,0);
}
\\
      };
\end{tikzpicture}
}
\quad
\resizebox{2.5cm}{!}{
\begin{tikzpicture}[baseline=-0.5ex]
      \tikzset{every left delimiter/.style={xshift=-1ex},every right delimiter/.style={xshift=1ex}}
      \matrix (mat) [left delimiter=(, right delimiter=)] {
      \draw;
\foreach \j in {0,4,...,28}{
    \draw[fill=gray!25,line width=10pt] (\j,32) rectangle (\j+4,0);
}
\draw[fill=gray!25,line width=10pt] (0,16) rectangle (32,16);
\\
      };
\end{tikzpicture}
}
\quad
\resizebox{2.5cm}{!}{
\begin{tikzpicture}[baseline=-0.5ex]
      \tikzset{every left delimiter/.style={xshift=-1ex},every right delimiter/.style={xshift=1ex}}
      \matrix (mat) [left delimiter=(, right delimiter=)] {
      \draw;
\foreach \j in {0,8,...,24}{
    \draw[fill=gray!25,line width=10pt] (\j,32) rectangle (\j+8,0);
}
\draw[fill=gray!25,line width=10pt] (0,8) rectangle (32,8);
\draw[fill=gray!25,line width=10pt] (0,16) rectangle (32,16);
\draw[fill=gray!25,line width=10pt] (0,24) rectangle (32,24);
\\
      };
\end{tikzpicture}
}
\quad
\resizebox{2.5cm}{!}{
\begin{tikzpicture}[baseline=-0.5ex]
      \tikzset{every left delimiter/.style={xshift=-1ex},every right delimiter/.style={xshift=1ex}}
      \matrix (mat) [left delimiter=(, right delimiter=)] {
      \draw;
\foreach \j in {0,4,...,28}{
    \draw[fill=gray!25,line width=10pt] (32,\j) rectangle (0,\j+4);
}
\draw[fill=gray!25,line width=10pt] (16,0) rectangle (16,32);
\\
      };
\end{tikzpicture}
}
\quad
\resizebox{2.5cm}{!}{
\begin{tikzpicture}[baseline=-0.5ex]
      \tikzset{every left delimiter/.style={xshift=-1ex},every right delimiter/.style={xshift=1ex}}
      \matrix (mat) [left delimiter=(, right delimiter=)] {
      \draw;
\foreach \j in {0,2,...,30}{
    \draw[fill=gray!25,line width=10pt] (32,\j) rectangle (0,\j+2);
}
\\
      };
\end{tikzpicture}
}
\\
\end{minipage}
\caption{Hierarchical decomposition of the row and column indices
    of a $16\times 16$ matrix. The dyadic trees $T_{X}$ and $T_{\Omega}$ have roots containing $16$ rows and $16$ columns respectively, and their leaves containing only a single row or column. The partition above indicates the complementary low-rank property of the matrix, and assumes that each submatrix is rank-$1$. }
\label{fig:complement}
\end{figure}

We first recall the definition of a complementary low-rank matrix given in \cite{BF}. 
Suppose that $K \in \bbC^{N \times N}$.   We denote the set of rows of $K$ by $X$ and the set of columns of $K$ by $\Omega$.
We introduce two trees $T_X$ and $T_\Omega$ that are generated by bisecting the sets $X$ and $\Omega$ recursively, and the elements of $T_X$ and $T_\Omega$ consist of subsets of $X$ and $\Omega$, respectively. Assume that both trees have the same depth $L=\O{\L{N}}$ with the top-level being level $0$ and the bottom one being level $L$ (see Figure~\ref{fig:domain-tree-BF} for an illustration). The top level of $T_X$ (and $T_\Omega$) contains all the indices in $X$ (and $\Omega$), while each leaf at the bottom level contains $O(1)$ indices. Then, the matrix $K$ is said to satisfy the {\it complementary low-rank property}, if the following property holds: for any level $\ell$, any node $A \in T_X$ at level $\ell$, and any node $B \in T_\Omega$ at level $L-\ell$, the submatrix $K(A,B)$, obtained by restricting $K$ to the rows indexed by the points in $A$ and the columns indexed by the points in $B$, is numerically low-rank.
By numerically low-rank, we mean that the ranks of the submatrices grow no more quickly than $\log^\kappa\left(N\right)$ with the size of the matrix $K$.
In many cases of interest, $\kappa=0$ --- that is,
the ranks of the submatrices are bounded by a constant independent of $N$.
See Figure~\ref{fig:complement} for an illustration of the complementary low-rank property.

An interpolative decomposition butterfly factorization (IDBF) of a complementary low-rank matrix $K$
is a factorization of the form
\begin{equation}
    K \approx U^{L}U^{L-1}\cdots U^{h} S^h V^{h}\cdots V^{L-1}V^{L},
\end{equation}
where $L$ is the number of levels in the trees $T_X$ and $T_\Omega$, $h = L/2$
and each of the matrices $U^l$ and $V^l$ is sparse with 
$\O{N}$ entries.  The number of levels $L$ in this decomposition is on the order  
of $\L{N}$, where $N$ is the dimension of $K$.
This factorization is obtained by constructing
interpolation decompositions of the low rank blocks of $K$ using the algorithm
of the preceding section.    The IDBF algorithm takes as input the 
parameter $r_k$ which is an estimate of the maximum possible ranks of the low-rank blocks.

Given an $N \times N$ matrix $K$, or equivalently an $\O{1}$ algorithm to evaluate an arbitrary entry of $K$, 
the algorithm of \cite{IDBF} constructs this data-sparse representation of $K$ 
in $\O{N \log^{\kappa+1}\left(N\right)}$ operations using 
$\O{N \log^{\kappa+1}\left(N\right)}$ storage.  
Once this factorization has been constructed,
the matrix $K$ can be applied in $\O{N \log^{\kappa+1} \left(N\right)}$  operations.

\subsection{Low-rank approximation by randomized sampling} \label{sec:rSVD}

In this section, we discuss an existing linear complexity algorithm for constructing
an approximate singular value decomposition (SVD) of a matrix.

Suppose that $A \in \bbC^{m\times n}$ has singular values
\begin{equation}
\left|\sigma_1\right| \geq \left|\sigma_2\right| \geq \cdots \geq \left|\sigma_l\right|,
\end{equation}
where $l = \min(n,m)$.  A rank-$r$ singular value decomposition of $A$ is a factorization
of the form
\begin{equation}\label{eqn:SVD2}
    A \approx U \Sigma V^T,
\end{equation}
where $U\in \bbC^{m\times r}$ is orthogonal, $\Sigma\in \bbR^{r\times r}$ is diagonal, 
 $V\in \bbC^{n\times r}$ is orthogonal and
\begin{equation}
\left\|A-U\Sigma V^T\right\|_2 = \sigma_{k+1}.
\end{equation}
The construction of a factorization of this form is a notoriously expensive calculation.  However,
using randomized algorithms, approximate SVDs of the same form with slightly lower accuracy
can be rapidly constructed.  That is, randomized algorithms result in factorizations of the 
form (\ref{eqn:SVD2}) for which 
\begin{equation}
\left\|A-U\Sigma V^T\right\|_2 
\end{equation}
is no longer equal to the optimal value $\sigma_{k+1}$, but is instead slightly larger.

One of the first practical randomized algorithms for constructing approximate SVDs was proposed in \cite{randomSVD}.  It operates by applying a random transform to the matrix $A$
and requires $\O{nmk}$ operations.  By using a special random transform
which can be applied rapidly via the FFT, 
a variant of the algorithm of \cite{randomSVD} which requires 
$\O{nm\L{k}}$ can be obtained.

In \cite{RandomSample}, a method which operates by 
 randomly sampling $\O{1}$  rows and columns of the input matrix
is described.  It  only requires $\O{m+n}$ operations and $\O{m+n}$ storage.
Here, we denote this algorithm as Function \rsvd and it is presented in Algorithm~\ref{alg:rSVD}. 
Assuming the whole low-rank matrix $A$ is known, the input of Function \rsvd is $A$, $\O{1}$ 
randomly sampled row indices $\mathcal{R}$ and column indices $\mathcal{C}$, as well the 
parameter $r$.  Equivalently, it can also be assumed that $A(\mathcal{R},:)$ 
and $A(:,\mathcal{C})$ are known as the inputs. 
The outputs are the matrices $U\in \bbC^{m \times r}$, $\Sigma\in \bbR^{r \times r}$, and $V\in \bbC^{n\times r}$ 
which give an approximate SVD \eqref{eqn:SVD2}. 

In Function \rsvd, for simplicity, given any matrix $K \in \bbC^{s\times t}$, 
Function \qr{K} performs a pivoted QR decomposition $K(:,P) = QR$, where $P$ is a permutation vector of the $t$ 
columns, $Q$ is a unitary matrix, and $R$ is an upper triangular matrix with positive diagonal entries in decreasing order. 
Function \rp{m,r} denotes an algorithm that randomly selects $r$ different samples in the set $\{1,2,\dots,m\}$. 

In most cases, to obtain higher accuracy, we add an oversampling parameter $q$ and we sample $rq$ rows and columns and only generate a rank $r$ truncated SVD in the penultimate line in Algorithm~\ref{alg:rSVD}. Larger $q$ results in 
better stability of Algorithm~\ref{alg:rSVD}. In our numerical experiments, Algorithm 1 is stable with high probability for $q >= 2$, and $q = 2$ is empirically sufficient to achieve accurate low-rank approximations.

\begin{algorithm2e}
    
    \Fn{$\left[U, \Sigma, V\right] \leftarrow$ \rsvd{$A, \mathcal{R}, \mathcal{C}, r$}}{
    
    $\left[m,n\right] \leftarrow \sz(A)$
    
    $P \leftarrow $ \qr{$A(\mathcal{R},:)$} \qd $\Pi_{col} \leftarrow P(1:r)$
    \tcp*[f]{$A(\mathcal{R},P) = QR$}
    
    $P \leftarrow $ \qr{$A(:,\mathcal{C})^T$} \qd $\Pi_{row} \leftarrow P(1:r)$
    \tcp*[f]{$A(P,\mathcal{C}) = R^T Q^T$}
    
    $Q \leftarrow $ \qr{$A(:,\Pi_{col})$} \qd $Q_{col} \leftarrow Q(:,1:r)$
    \tcp*[f]{$A(P,\Pi_{col}) = QR$}
    
    $Q \leftarrow $ \qr{$A(\Pi_{row},:)^T$} \qd $Q_{row} \leftarrow Q(:,1:r)$
    \tcp*[f]{$A(\Pi_{row},P) = R^T Q^T$}
    
    $S_{row} \leftarrow $ \rp{$m,r$} \qd $I \leftarrow [\Pi_{row}, S_{row}]$
    
    $S_{col} \leftarrow $ \rp{$n,r$} \qd $J \leftarrow [\Pi_{col}, S_{col}]$
    
    $M \leftarrow \left( Q_{col}(I,:) \right) ^\dagger A(I,J) \left( Q_{row}^T(:,J) \right) ^\dagger$
    \tcp*[f]{$(\cdot)^\dagger:$ pseudo-inverse}
    
    $\left[U_M, \Sigma_M, V_M \right]\leftarrow \svd(M)$ \label{alg:rSVD:svd}
    
    $U \leftarrow Q_{col} U_M$ \qd $\Sigma \leftarrow \Sigma_M$ \qd $V \leftarrow Q_{row} V_M$
    }
    
    \caption{Randomized sampling for a rank-$r$ approximate SVD with $\O{m+n}$ operations, such that $A \approx U \Sigma V^T$.}
    \label{alg:rSVD}
    
\end{algorithm2e}

\section{Algorithm for the application of the ALT} \label{sec:algo}

The principal content of this section is a description of our 
block partitioning algorithm based on IDBF and low-rank approximation by randomized 
sampling for applying the forward ALT.
Before we present this, however, we 
briefly discuss certain background information regarding the associated Legendre functions and the
associated Legendre transform which we exploit.
As we observe in Section~\ref{sec:inverse}, the inverse ALT can be applied in essentially the same fashion
as the forward ALT.

\subsection{Background}

\subsubsection{The relationship between the forward and inverse associated Legendre transforms}
\label{sec:inverse}

For fixed $N$ and $|m|\leq N$, the forward ALT consists of computing the values of the
sum
\begin{equation}
\label{eqn:new1}
g(m,\theta) = \sum_{k=|m|}^{2N-1} \beta_{k,m} \overline{P}_k^{|m|} (\cos(\theta)),
\end{equation}
at the nodes of the $2N$-point Gauss-Legendre quadrature rule.    
We let 
\begin{equation*}
x_0 = \cos\left(\theta_0\right), x_1 = \cos\left(\theta_1\right), \ldots, x_{2N-1}=\cos\left(\theta_{2N-1}\right)
\end{equation*}
and
\begin{equation*}
w_0, w_1,\ldots, w_{2N-1}
\end{equation*}
denote the nodes and weights of this quadrature.  There are $(2N-|m|)$ coefficients in the expansion
(\ref{eqn:new1}) and $2N$ target points, so this amounts to applying the $2N \times (2N-|m|)$ matrix
whose $ij$ entry is 
\begin{equation*}
\overline{P}_j^{|m|} (\cos(\theta_i))
\end{equation*}
to the vector
\begin{equation}
\left(
\begin{array}{c}
\beta_{|m|,m} \\
\beta_{|m|+1,m} \\
\vdots\\
\beta_{2N-1,m} \\
\end{array}
\right).
\label{eqn:new2}
\end{equation}
of coefficients.

It is well-known that for $k \geq |m|$,  $\overline{P}_k^{|m|} (x)$ is a polynomial of degree $k-|m|$, and that
the functions 
\begin{equation}
\left\{ \overline{P}_k^{|m|} (x) : k=|m|, \ldots, 2N-1 \right\}
\end{equation}
form an orthonormal basis in the space of polynomials of degree no larger than $2N-1$.
The $2N$-point Gauss-Legendre quadrature rule exactly integrates the product of any two polynomials of degree $2N-1$.   In particular, it follows that  
when the $(2N-|m|) \times 2N$ matrix whose $ij$  entry is 
\begin{equation*}
\overline{P}_i^{|m|} (\cos(\theta_j)) w_j
\end{equation*}
is applied to the vector
\begin{equation*}
\left(
\begin{array}{c}
g(m,\theta_0)\\
g(m,\theta_1)\\
\vdots \\
g(m,\theta_{2N-1}) \\
\end{array}
\right)
\end{equation*}
the result is the vector of coefficients (\ref{eqn:new2}).
In other words, due to the orthonormality of the associated Legendre polynomials
and the method used to discretize spherical harmonic transforms, the matrix $B$ which
discretizes the inverse ALT is related to the matrix $A$ discretizing the forward ALT via
the formula
\begin{equation}
B = A^T W,
\end{equation}
where  $W$ is a diagonal matrix.    The methodology described in this section
for applying the forward ALT can be easily used to apply its transpose,
and hence also the inverse ALT.

\subsubsection
{Odd and even Legendre transform matrices} \label{sec:oddeven}

It is well-known (see, for instance, Chapter~14 of \cite{DLMF}) that 
$ \overline{P}_k^m(x)$ is odd when $k-|m|$ is odd, and even when
 $k-|m|$ is even.  This, together with the fact that the Gauss-Legendre
quadrature nodes are symmetric around $0$, allows us to 
reduce the cost of applying the forward ALT by a factor of $2$. 

More explicitly, the sum \eqref{eqn:new1} can be rewritten as 
\begin{equation}
    \label{eqn:sht4}
        g(\theta,m) = g_1(\theta,m) + g_2(\theta,m),
\end{equation}
where $g_1$ and $g_2$ are defined via the formulas
\begin{equation}
g_1(\theta,n)    = \sum_{0 \le k \le 2N-|m|-1,\, k \text{ is odd}} \beta_{k+|m|,m} \overline{P}_{k+|m|}^{|m|} (\cos(\theta)) 
\label{eqn:new3}
\end{equation}
and
\begin{equation}
g_2(\theta,n) = \sum_{0 \le k \le 2N-|m|-1,\, k \text{ is even}} \beta_{k+|m|,m} \overline{P}_{k+|m|}^{|m|} (\cos(\theta)).
\label{eqn:new4}
\end{equation}
Because of the symmetry of the Gauss-Legendre nodes, we have 
\begin{equation}
g_1(\theta_l,m) = -g_1(\theta_{2N-1-l},m) \ \ \mbox{and} \ \ \ g_2(\theta_l,m) = g_2(\theta_{2N-1-l},m)
\end{equation}
for  $l = 0, 1, \dots, 2N-1$.   Therefore, we can reduce the cost of applying the forward ALT by only computing
the values of \eqref{eqn:new1} and \eqref{eqn:new2} at the nodes $\theta_0,\theta_1,\ldots,\theta_{N-1}$
and using these to compute $g(m,\theta)$ at each of the Gauss-Legendre nodes.

Computing the sum (\ref{eqn:new3}) at each of the $N$ positive Gauss-Legendre nodes amounts to applying
an  $N \times \left(N - \ceil{\frac{|m|}{2}} \right)$ matrix, which we refer to as the 
odd ALT matrix.  Computing \eqref{eqn:new3} at each of the $N$ positive Gauss-Legendre nodes amounts
to applying  an $N \times \left(N - \floor{\frac{|m|}{2}} \right)$, which we refer to as the even ALT matrix.

\subsection{A block partitioning scheme} \label{sec:osci}

When $|m| >0$, the associated Legendre function $\overline{P}_k^{|m|}(\cos(\theta))$ has a single turning point (or the first inflection point) on the interval $(0,\pi/2)$.  Its location is given by the formula
\begin{equation}
    t^*_{k,m} = \arcsin{\left(\frac{\sqrt{m^2-1/4}}{k+1/2}\right)}.
\label{eqn:new5}
\end{equation}
See, for instance, Chapter~14 of \cite{DLMF} for details.  On the interval $(0,t^*_{k,m})$, 
$\overline{P}_k^{|m|}(\cos(\theta))$ is nonoscillatory and on $(t^*_{k,m},\pi/2)$
it is oscillatory.
We can view \eqref{eqn:new5} as defining a piecewise continuous curve that divides
the odd and even ALT matrices into oscillatory and nonoscillatory regions.  We will refer to this
as the ``turning point curve''.
Any subblock of these matrices that intersects this curve will have a high rank.
As a consequence of this, the even and odd ALT matrices do not have the complementary low-rank property.
Figure~\ref{fig:lg}~(a) shows an example of an odd ALT matrix with a graph of this piecewise continuous curve overlaid on top of it.

We use the following procedure to hierarchically partition the even and odd and ALT matrices
into blocks each of which will either be of small dimension or will have complementary low-rank property.
We will denote the resulting collection of subblocks by $\mathcal{B}_s$.
Since the shape of the matrices are not square when $m \neq 0$, we initially take $\mathcal{B}_s$
to consist of blocks each of which consist of all columns of the matrix and $b$ rows, where
\begin{equation}
    b = \left\{
    \begin{array}{ll}
        \lfloor \frac{N}{N - \ceil{\frac{|m|}{2}}} \rceil & \text{for odd matrices}, \\
        \lfloor \frac{N}{N - \floor{\frac{|m|}{2}}} \rceil & \text{for even matrices}.
    \end{array}
    \right.
    \label{eqn:nb}
\end{equation}
Each of these blocks is nearly square.
The symbol $\lfloor x \rceil$ means the nearest integer to $x$. 
Next, we repeatedly apply the following procedure.
We split each block in $\mathcal{B}_s$ which intersects the piecewise curve defined by \eqref{eqn:new5}
into a $2 \times 2$ grid of subblocks.  We stop this procedure only once 
all blocks which contain turning points have either fewer than $n_0$ rows or columns, where $n_0$
is a specified parameter.
This makes  the  maximum partition level is $L = \L{\frac{N}{n_0}}$. 
For each partition level $\ell$ of $\mathcal{B}_s$, the turning point curve intersects no 
more than $2^\ell - 1$ submatrices. 
Therefore, this procedure takes at most $\O{N}$ operations to partition the odd and the even matrices 
into submatrices, since
\begin{equation}
    \O{\sum_{\ell=1}^L \left( 2^\ell - 1 \right)} \sim \O{\frac{2N}{n_0} - \L{\frac{N}{n_0}} - 2} \sim \O{N}.
\end{equation}
At level $\ell$, the sub-matrix is of size $\O{\frac{N}{2^\ell}}\times \O{\frac{N}{2^\ell}}$ and there are $\O{2^\ell}$ such sub-matrices either as an oscillatory block or a non-oscillatory block. See Figure~\ref{fig:lg}~(b), which shows an example of an odd ALT matrix that has been partitioned into blocks by this
procedure. 

Finally, we estimate the cost for applying matrix blocks that intersect with the turning point curve.
The maximum partition level is $L = \L{\frac{N}{n_0}}$ and the turning point curve
intersects with no more than $2^\ell - 1$ submatrices for each partition level $\ell$.
Each box has size at most $n_0 \times n_0$ implying that the cost of applying all these blocks is 
\begin{equation}
\label{eqn:num_blocks}
\O{n_0^2 L \left(2^{L}-1\right)}
= 
\O{N \L{N}}.
\end{equation}

\subsection{Factorization and application of matrix blocks} \label{sec:fa}

In each of the partitioned matrices, there are three types of blocks:
oscillatory blocks,  non-oscillatory blocks, and the blocks which intersect the turning point curve. 
We deal with each of these different kinds of blocks by different approaches. 
In the following discussion, we assume that, for a fixed dimension $N$ of the ALT, the ranks of all low rank matrices are 
bounded.  We denote the  least upper bound by $r(N)$.
\begin{enumerate}
    \item For an oscillatory block $\mathcal{B}^o$ with the size $N_0 \times N_0$, we use the IDBF to construct a factorization
    \begin{equation}\label{eqn:IDBFf}
        \mathcal{B}^o \approx U^{L}U^{L-1}\cdots U^{h} S^h V^{h}\cdots V^{L-1}V^{L},
    \end{equation}
    where $L = \O{\L{N_0}}$ and $h = \frac{L}{2}$.  This takes only $\O{\left(r(N)\right)^2 N_0 \L{N_0}}$
operations and memory.  After
factorization, we can apply the subblock $B^o$ with $\O{(r(N))^2 N_0 \L{N_0}}$ operations and memory.
    
    \item 
In the non-oscillatory region, the entries of the odd and even ALT matrices can be of extremely
small magnitudes.  Therefore, when processing a non-oscillatory block $\mathcal{B}^n$ with the size $N_0 \times N_0$, we first take the largest subblock  $\mathcal{B}^{n'}$which does not contain
any elements of magnitude smaller than machine precision.
 Next, we use the algorithm of Section~\ref{sec:rSVD} to construct a low rank factorization of the form
    \begin{equation}\label{eqn:lrf}
        \mathcal{B}^{n'} \approx U \Sigma V^T
    \end{equation}
with $\O{(r(N))^2N_0 }$ operations  and $\O{ r(N) N_0 }$ memory.
After factorization, the application of  $\mathcal{B}^n$ requires $\O{ r(N) N_0 }$ operations  and memory.
    
    \item For a block $\mathcal{B}^t$ including turning points, we also let $\mathcal{B}^{t'}$ be a smaller submatrix
which excludes as many entries whose magnitudes are smaller than machine precision as possible.
These blocks are applied to a vector through a standard matrix-vector multiplication with no made attempt to accelerate it. The operation{s} and memory complexity for these blocks are bounded by the number of nonzero entries in these blocks, which $\O{N \L{N}}$ for an ALT matrix of size $N\times N$ according to \eqref{eqn:num_blocks}.
\end{enumerate}

\begin{figure}[!ht]
    \centering
    \begin{tabular}{ccc}
        \includegraphics[height=2.8in,trim={7.5cm 1cm 7cm 2cm},clip] {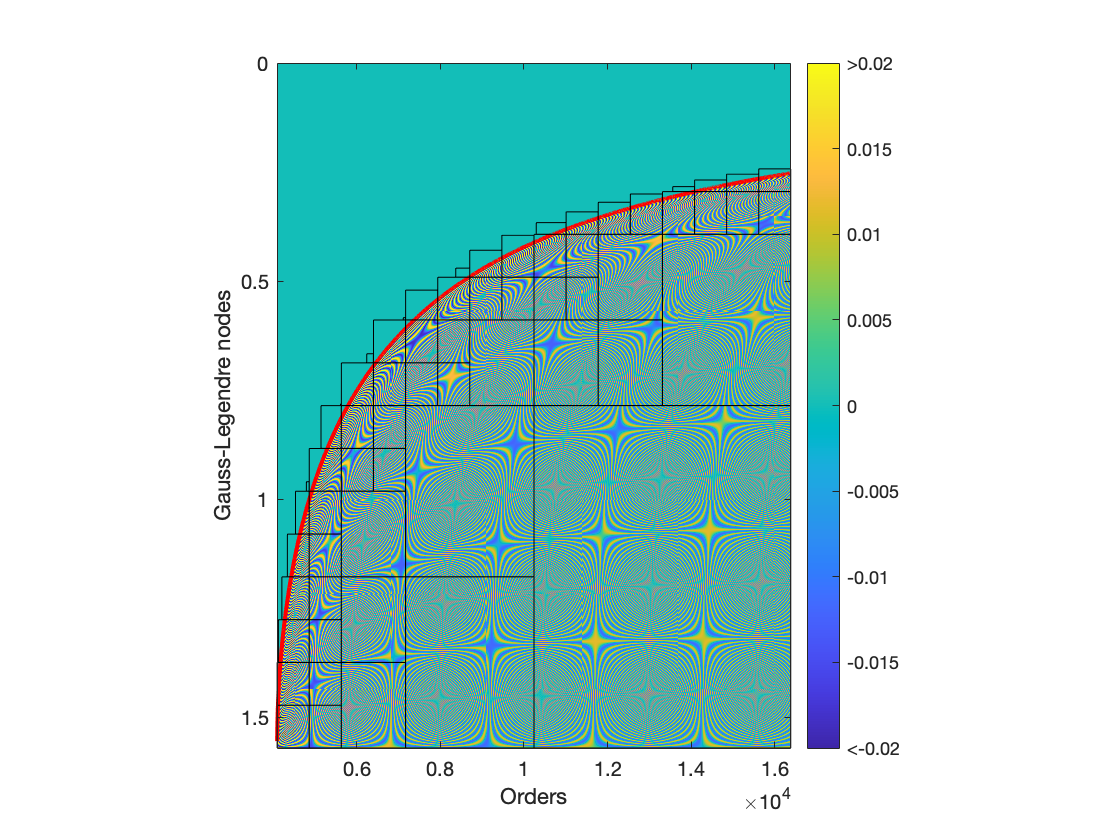} & \includegraphics[height=2.8in,trim={10.5cm 1cm 9.5cm 2cm},clip] {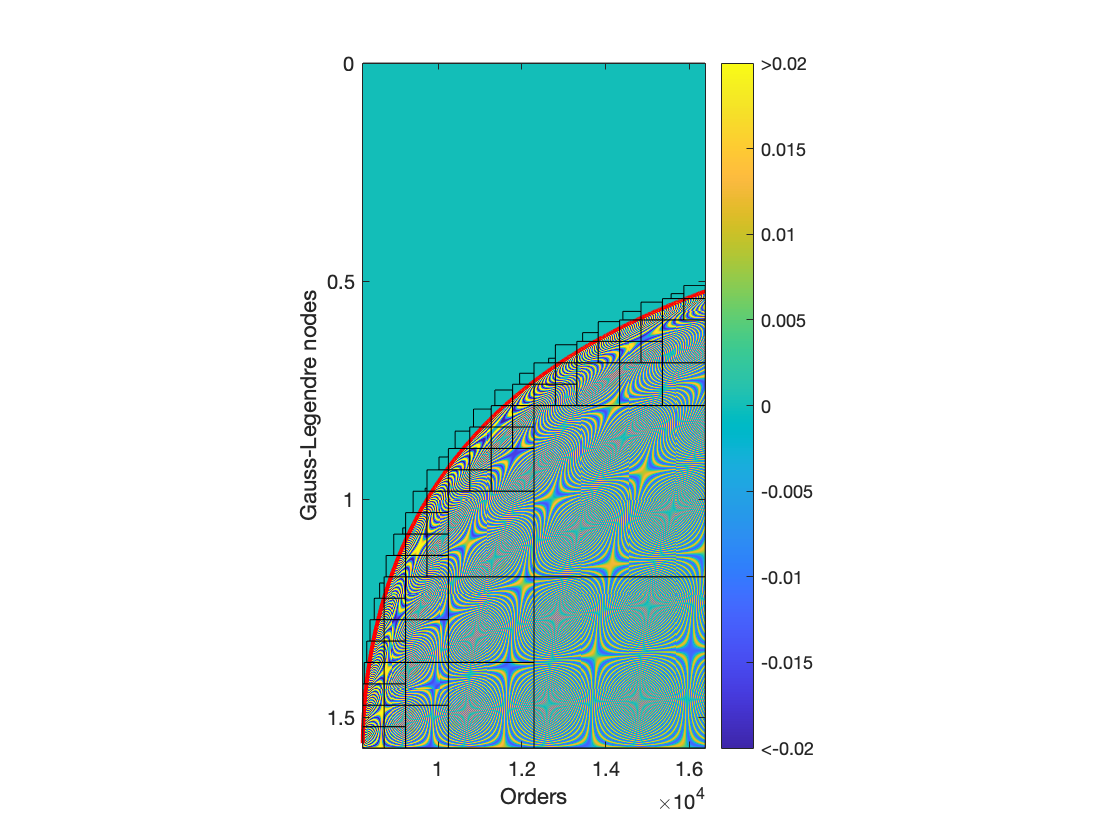} &  \includegraphics[height=2.8in,trim={13.5cm 1cm 14cm 2cm},clip] {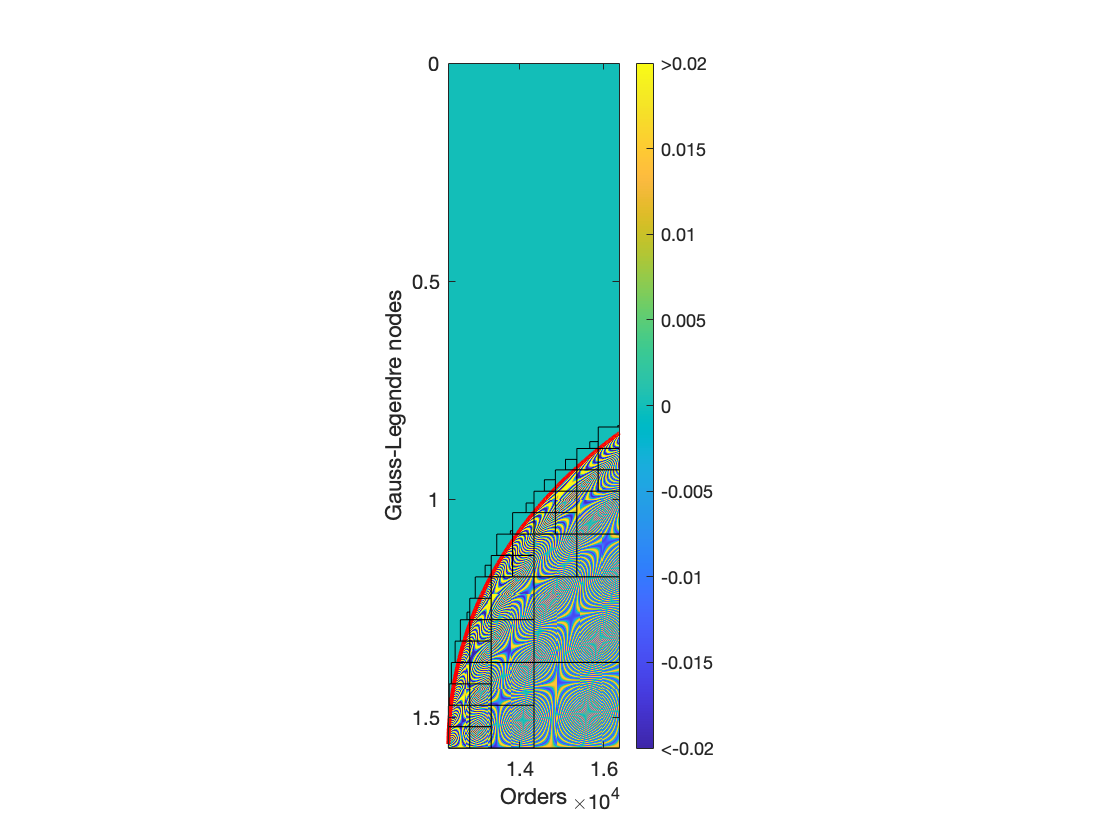}  \\
        (a) $m=4096$ &
        (b) $m=8192$ & (c) $m=12288$
    \end{tabular}
    \caption{A visualization of the partitioning procedure for the forward ALT matrix. 
In each case,  $N = 8192$. From left to right, the orders of the transform are
$m=4096$, $m=8192$ and $m=12288$.}
    \label{fig:8192}
\end{figure}

Figure~\ref{fig:8192} shows the boxes which result after as many elements of negligible magnitude as possible
have been excluded.  Empty blocks with $0 \times 0$ size are omitted in the figure and will not be utilized in the application step.

\section{Computational complexity}
\label{sec:analysis}

A rigorous estimate of the computational complexity of our algorithm would seem to require a 
bound on the ranks of the subblocks of the odd and even ALT matrices which are in the oscillatory
regions.  To the author's knowledge, no such bounds are presently known, except in the special case $m=0$ 
(such an estimate can be found in \cite{fastSHT}).
We can, however, develop an estimate on the complexity of our algorithm 
in terms of both operation{s} and memory 
assuming that the ranks of these boxes are bounded by a quantity depending on $N$, which
we denote by $r(N)$.

We will assume that the matrix we are applying is an $N \times N$ odd ALT matrix which we will
denote by $A$. The analysis for the even case is similar.
We first observe that subdividing a matrix with the complementary low-rank property into subblocks
and using the IDBF to apply each subblock separately has essentially the same asymptotic complexity
as using the IDBF to apply the entire matrix. Recall that there are $L = \L{\frac{N}{n_0}}$ levels of subdivision. 
At level $\ell$, the sub-matrix is of size $\O{\frac{N}{2^\ell}}\times \O{\frac{N}{2^\ell}}$ and 
there are $\O{2^\ell}$ such sub-matrices as an oscillatory block for IDBF. Hence, the total factorization 
and application complexity in terms of both operation{s} and memory is bounded by:
\[
\sum_{\ell=1}^{\L{\frac{N}{n_0}}} \O{2^\ell} ( (r(N))^2 \O{\frac{N}{2^\ell} 
\L{\frac{N}{2^\ell}}}=\O{(r(N))^2 N \Ltwo{N}}.
\] 

The complexity analysis is similar for non-oscillatory blocks. The numbers of blocks of different sizes are the same as those of oscillatory blocks. The factorization and application complexity in terms of operation{s} and memory of each non-oscillatory block is cheaper than that of an oscillatory block. Therefore, the total cost for non-oscillatory blocks has a scaling smaller than that of oscillatory blocks. In fact, due to the very rapid decay of the associated Legendre functions as one moves away from
the turning point into the nonoscillatory regime (see Figure~\ref{fig:8192}), many non-oscillatory blocks are almost zero and, hence, we neglect their computation.

Based on extensive numerical experiments, some of which are presented in the following section, we conjecture
that  $(r(N))^2$ grows as $\O{\L{N}}$.  Of course, if this is correct, then our algorithm
for applying the ALT requires $\O{N \Lthree{N}}$ operations and can be used to apply the SHT in 
$\O{N^2 \Lthree{N}}$ operations.

We summarize the situation with the following theorem and conjecture:

\begin{theorem}
\label{thm}
{For a fixed precision $\epsilon$}, the time and memory complexities of the  algorithm of Section~\ref{sec:algo} for applying
an odd or even ALT matrix of size $N\times N$ to a vector are 
$$
\O{(r(N))^2 N \Ltwo{N}},
$$
 where $r(N)$ is the least upper bound for both the ranks of the non-oscillatory blocks which do not 
intersect the turning point curve and the ranks of the subblocks of the oscillatory portion of the matrix.
\end{theorem}

\begin{conjecture}
\label{conj}
{For a fixed precision $\epsilon$}, the quantity $(r(N))^2$ grows as $\O{\L{N}}$.
\end{conjecture}

\section{Numerical results} \label{sec:results}

This section presents several numerical experiments which demonstrate the efficiency of the proposed algorithm. 
Our code was written in MATLAB\textsuperscript{\textregistered} and executed on a single 3.2GHz core.
It is available as a part of the ButterflyLab package (\url{htps://github.com/ButterflyLab/ButterflyLab}).

For the forward ALT, given an order $m$, we let $g^d(\theta)$ denote the results 
given by applying the discretized operator directly using a standard matrix-vector multiplication,
and we let 
 $g^b(\theta)$ 
denote the results obtained via the proposed algorithm.
The accuracy of our method is estimated via the relative error defined by
\begin{equation}
    \epsilon^{fwd} = \sqrt{\cfrac{\sum_{\theta\in S_1}|g^b(\theta)-g^d(\theta)|^2} {\sum_{\theta\in S_1}|g^d(\theta)|^2}},
\end{equation}
where $S_1$ is an index set containing $256$ randomly sampled row indices of the non-zero part in 
the odd matrix or the even matrix. 

We use $a^d(k)$ and $a^b(k)$ denote the results obtained by applying the inverse ALT using a standard
matrix-vector multiplication and via our algorithm, respectively.
The definition of the error $\epsilon^{inv}$ in this case is
\begin{equation}
    \epsilon^{inv} = \sqrt{\cfrac{\sum_{k\in S_2}|a^b(k)-a^d(k)|^2} {\sum_{k\in S_2}|a^d(k)|^2}},
\end{equation}
where $S_2$ is an index set containing $256$ randomly sampled row indices of the odd matrix or the even matrix.

In all of our examples, the tolerance parameter $\epsilon$ for interpolative decompositions is set to $10^{-10}$, 
the minimum length $n_0$ for the partitioned block is set to $512$, and the rank parameter $r$ 
for randomized SVD in low-rank phase matrix factorization is set to $30$. 

\begin{figure}[!ht]
    \centering
    \begin{tabular}{ccc}
        \includegraphics[height=2.1in,trim={3.5cm 6.3cm 4.2cm 7cm}, clip]{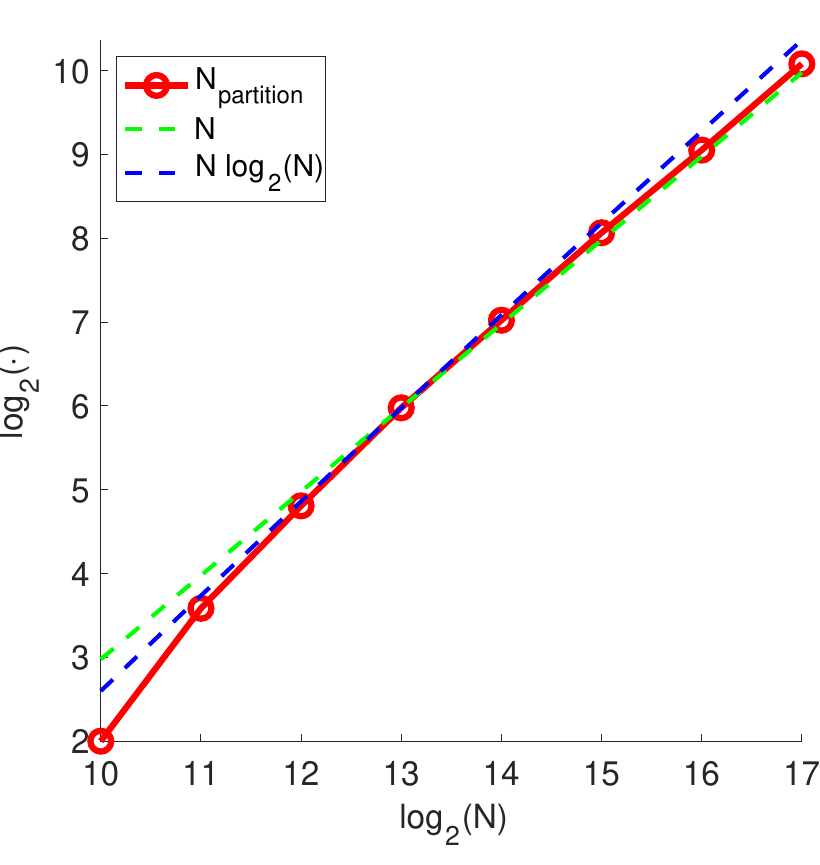} &
        \includegraphics[height=2.1in,trim={3.5cm 6.3cm 4.2cm 7cm}, clip]{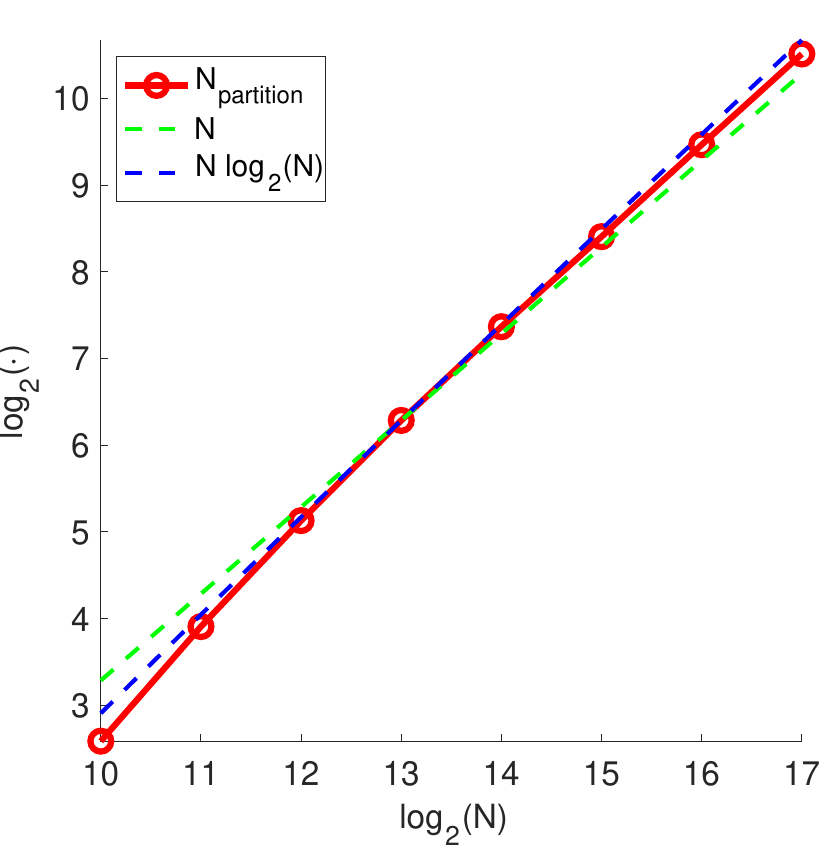} &
        \includegraphics[height=2.1in,trim={3.5cm 6.3cm 4.2cm 7cm}, clip]{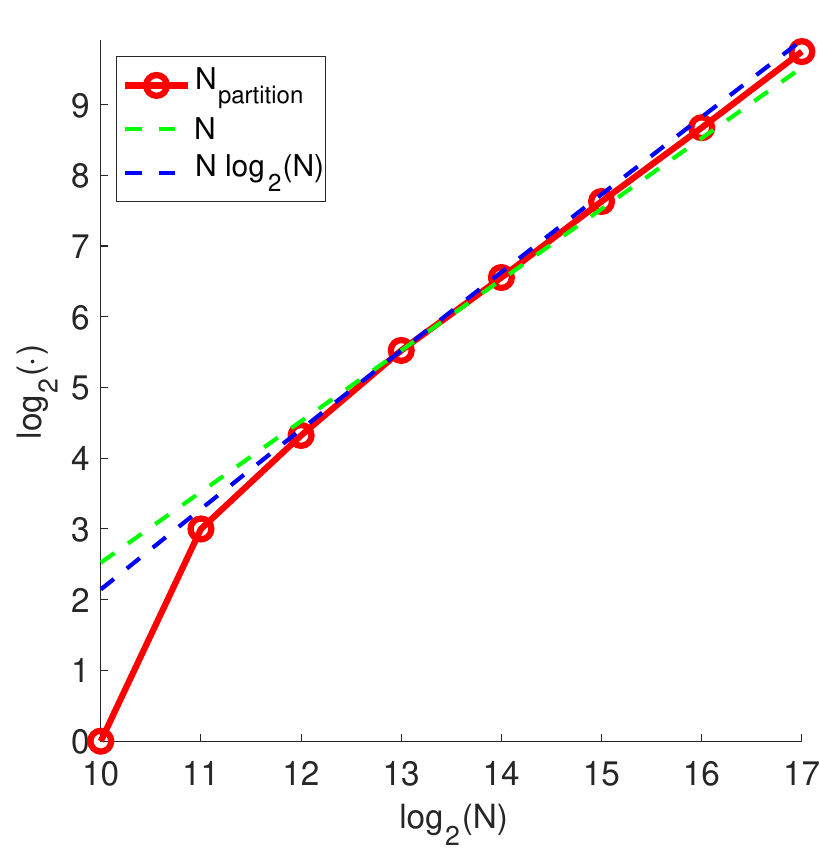} \\
        (a) $m = 0.5N$ & (b) $m = N$ & (c) $m = 1.5N$
    \end{tabular}
    \caption{Plots of the number of the remaining blocks a a function of $N$ for $m = 0.5N$, $m = N$, and $1.5N$. {The reference functions $N$ and $N \L{N}$ are multiplied by a constant in order to fit the data.}}
    \label{fig:nnp}
\end{figure}

\paragraph{Number of the blocks:}
Our first experiment consists of counting the number of blocks which remain after those which contain
no non-negligible elements are discarded.
Figure~\ref{fig:nnp} visualizes the results of this experiment for different $N$ and with 
$m$ set to be $0.5N$, $N$ and $1.5N$.  We observe that the number of remaining blocks scales nearly linearly as the problem size increases.

\begin{figure}[!ht]
    \centering
    \begin{tabular}{c}
        \includegraphics[width=11cm,trim={1cm 5.5cm 2.5cm 6cm}, clip]{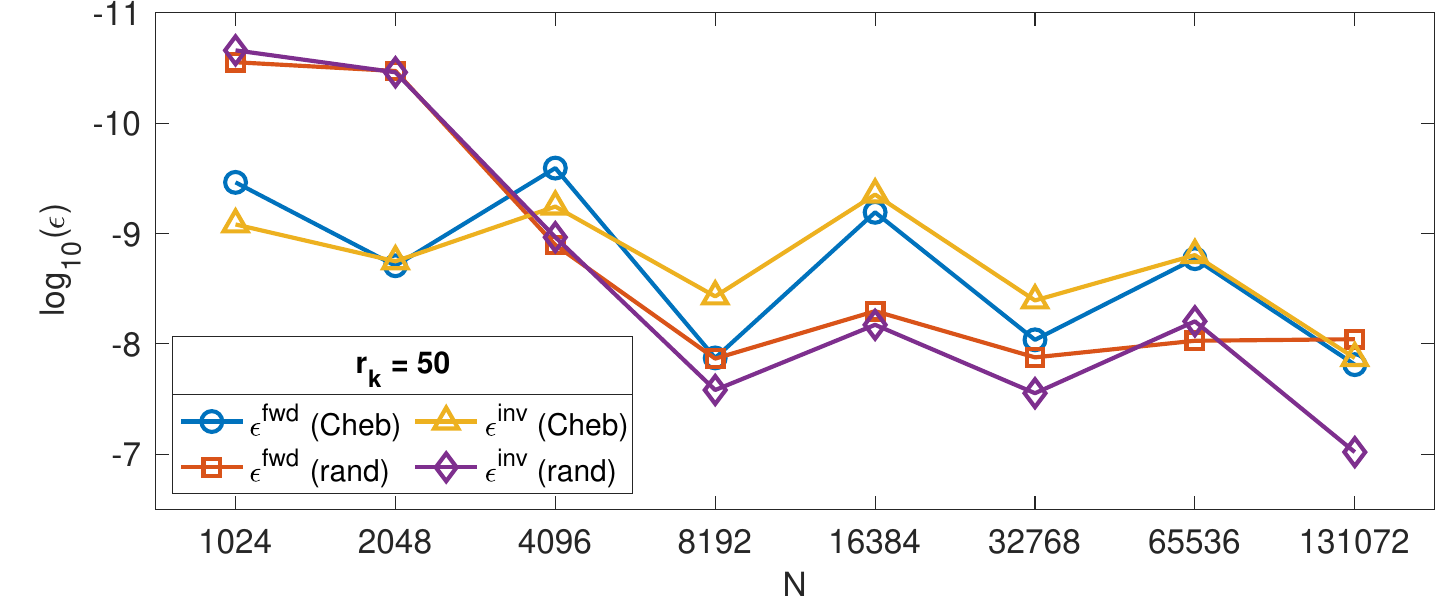} \\
        \includegraphics[width=11cm,trim={1cm 5.5cm 2.5cm 6cm}, clip]{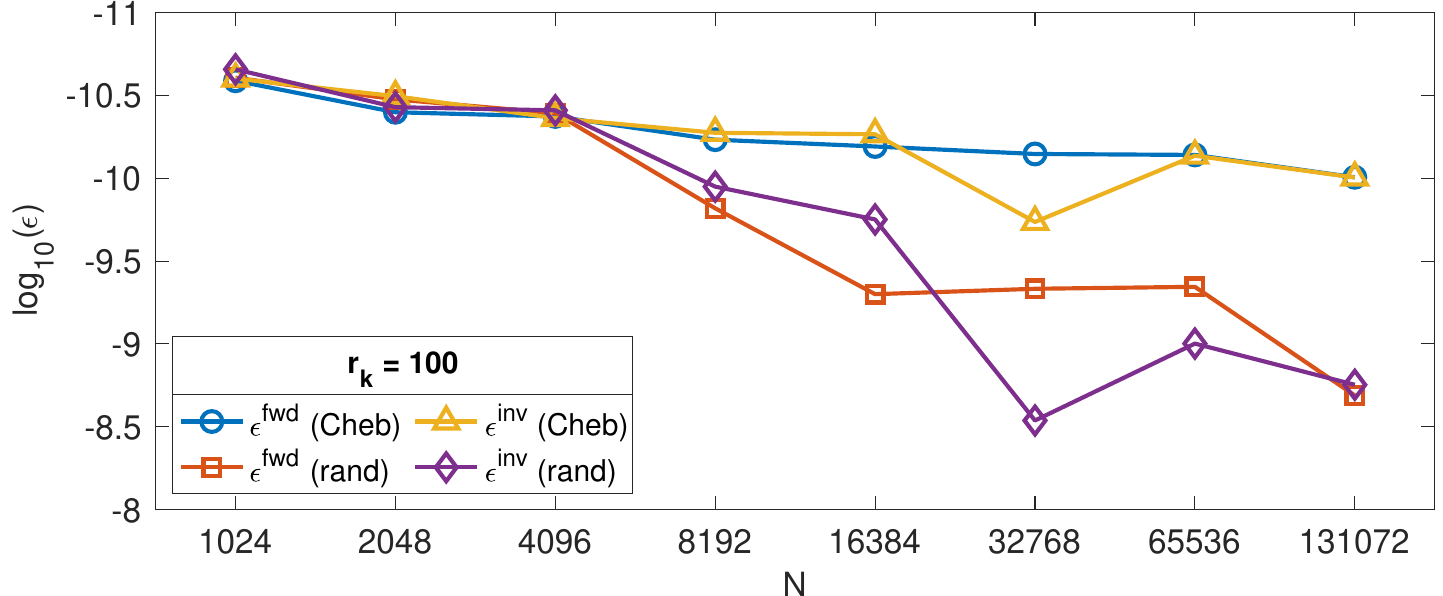} \\
        \includegraphics[width=11cm,trim={1cm 5.5cm 2.5cm 6cm}, clip]{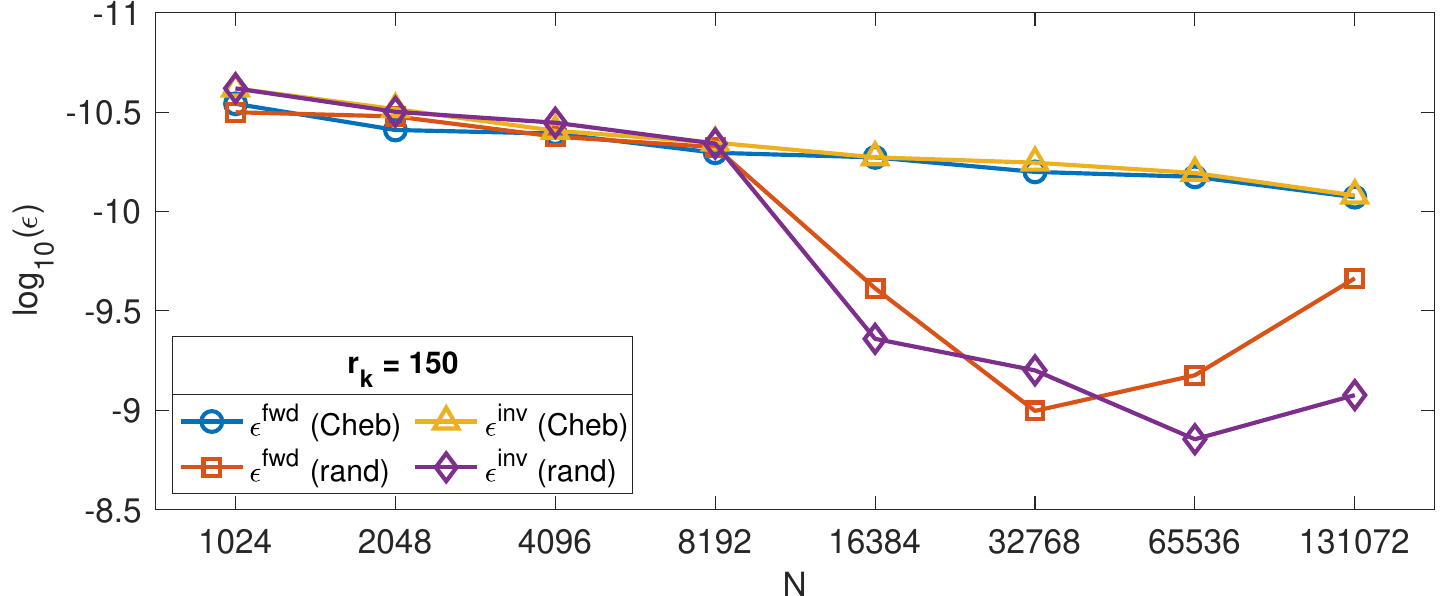}
    \end{tabular}
    \caption{The results of experiments comparing the error in applying the associated Legendre transform 
when different grids of points are used to form interpolative decompositions in the IDBF algorithm.
Here, $N$ is the size of the matrix, and the order $m$ is set to be $N$ in each case. 
The adaptive rank $r_k$ for IDBF is set to be $50$, $100$ and $150$ from the top panel to the bottom panel. 
``Cheb" and ``rand" represent IDs with Mock-Chebyshev points and randomly selected points, respectively.}
    \label{fig:point}
\end{figure}

\paragraph{Selection of Mock-Chebyshev points or randomly selected points:}
Next, we compare the results of using  Mock-Chebyshev points and randomly selected points for evaluating 
IDs in the IDBF process.
The results are shown in Figure~\ref{fig:point}.  In these experiments, the order parameter $m$ is set 
to be equal to $N$ and the adaptive rank $r_k$ for IDBF is set to be $50$, $100$ or $150$.
We observe that the accuracy of results increases as the rank parameter $r_k$ increases,
and that the accuracy of IDs performed with Mock-Chebyshev points is higher than that of the IDs performed
with randomly selected points.  Moreover, we conclude that letting $r_k = 150$ 
suffices to achieve high-accuracy with Mock-Chebyshev points. 

Thus, for the rest of 
experiments, we will use Mock-Chebyshev points as grids to compute IDs in the IDBF algorithm, 
and the adaptive rank $r_k$ for IDBF will be fixed at $150$.

\paragraph{Associated Legendre transforms of different orders:}
In these experiments, we measured the accuracy and efficiency of the proposed algorithm for various orders of $m$.

\begin{figure}[!ht]
    \centering
    \begin{tabular}{c}
        \includegraphics[width=11cm,trim={1cm 4.3cm 2.5cm 5.1cm}, clip]{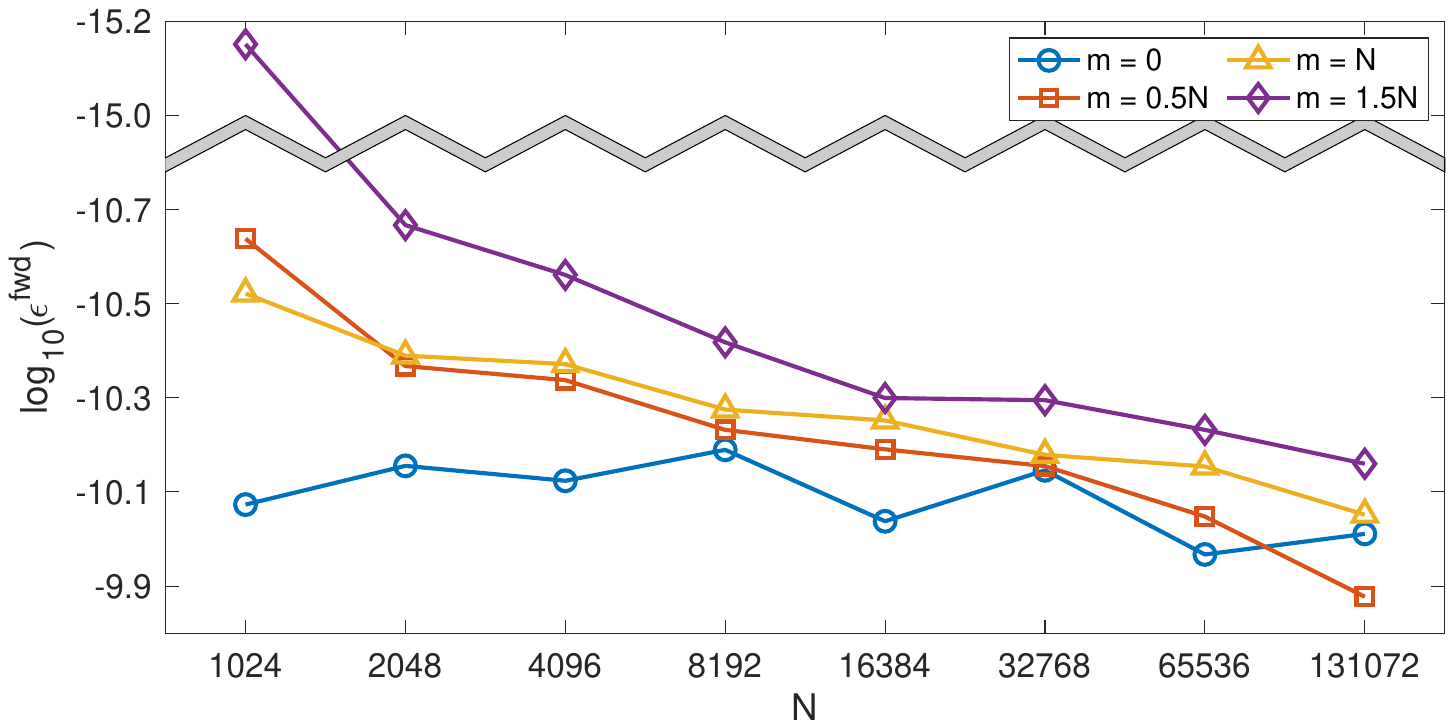} \\
        \includegraphics[width=11cm,trim={1cm 4.3cm 2.5cm 5.1cm}, clip]{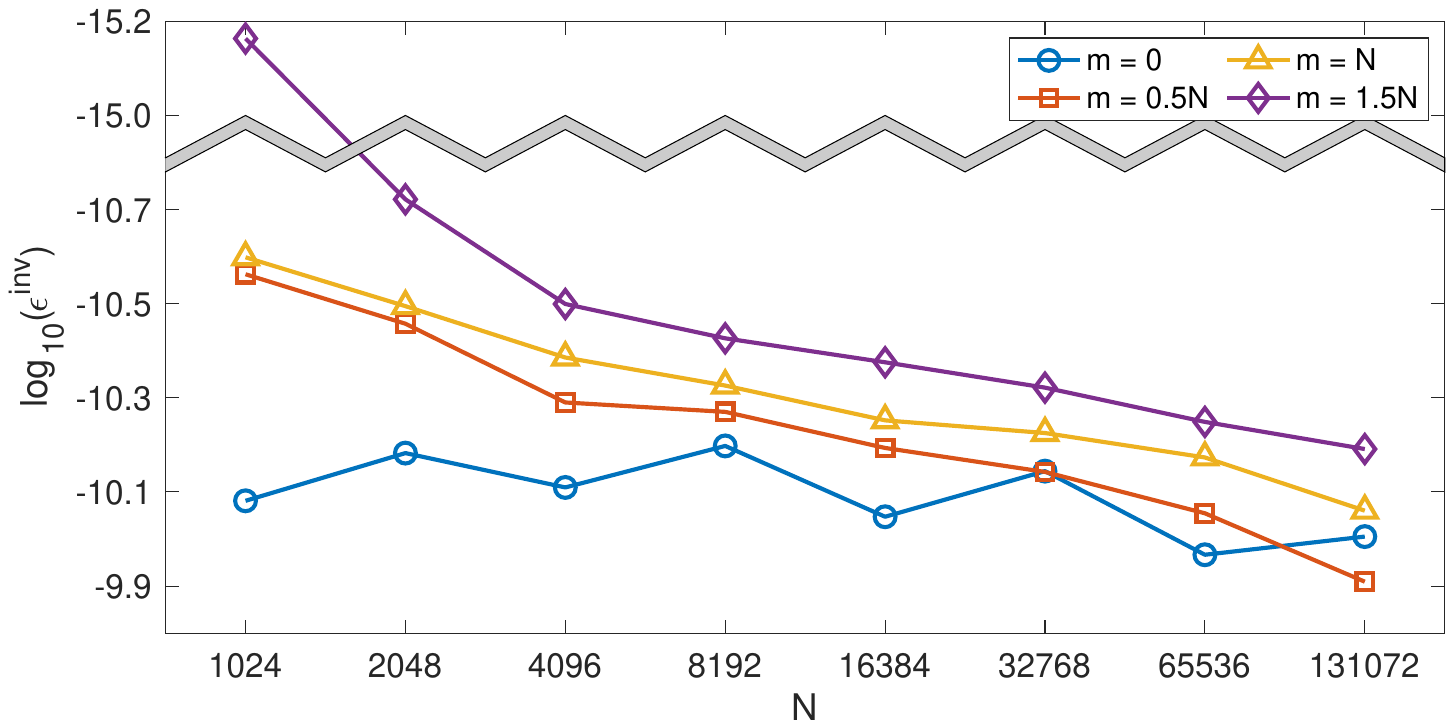}
    \end{tabular}
    \caption{The  errors in the application of associated  Legendre transforms
of  different orders.   The top panel shows the error $\epsilon^{fwd}$ of the forward ALT, 
and the bottom panel shows that error $\epsilon^{inv}$ of the inverse ALT. Here, $N$ is the size of the matrix, 
and the order $m$ was set to be $0$, $0.5N$, $N$ or $1.5N$. The adaptive rank $r_k$ for IDBF was taken to be
 $150$.}
    \label{fig:error}
\end{figure}

Figure~\ref{fig:error} shows that the accuracy of the proposed algorithm is unaffected by
the order $m$ of the ALT, even though the accuracy decays slightly as the problem size increases. 
The slightly increasing error appears to be due to the randomness of the proposed algorithm in Subsection~\ref{sec:rSVD}.  As the problem size increases, the probability of capturing the low-rank 
matrix with a fixed rank parameter becomes slightly smaller.

\begin{figure}[!ht]
    \centering
    \begin{tabular}{cc}
        \multicolumn{2}{c}{\includegraphics[width=9cm,trim={3cm 16.2cm 8cm 4.2cm}, clip] {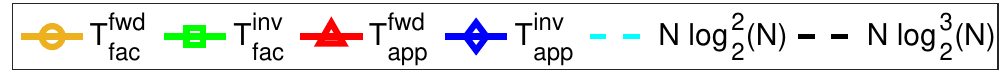}} \\
        \includegraphics[width=5cm,trim={7cm 3.3cm 7cm 5.39cm}, clip] {figure/00n/00n.pdf} &
        \includegraphics[width=5cm,trim={7cm 3.3cm 7cm 5.39cm}, clip] {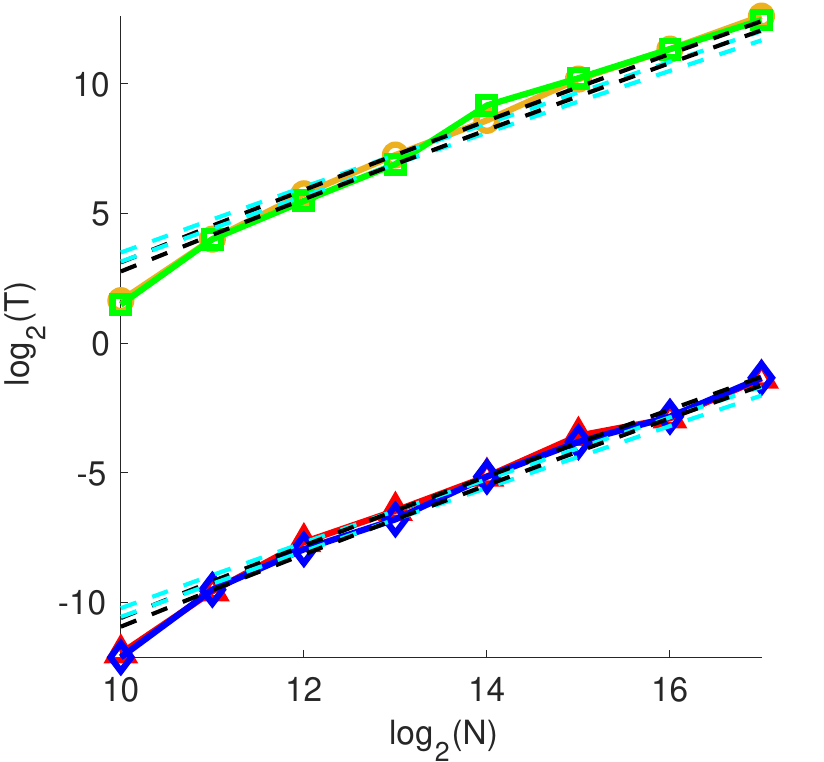} \\
        (1) $m=0$ & (2) $m=0.5N$ \\
        \includegraphics[width=5cm,trim={7cm 3.3cm 7cm 5.39cm}, clip] {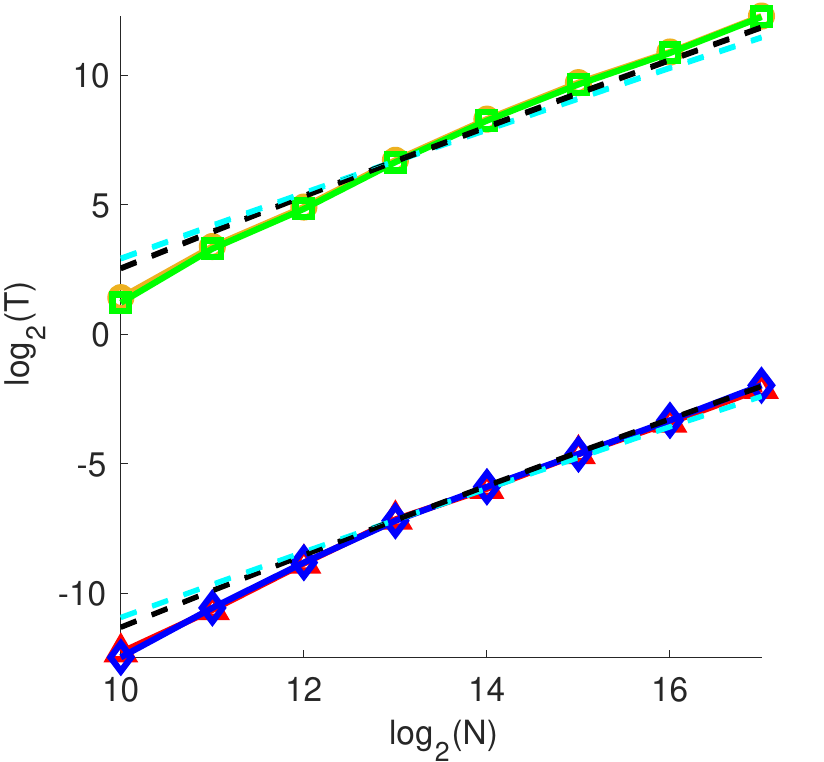} &
        \includegraphics[width=5cm,trim={7cm 3.3cm 7cm 5.39cm}, clip] {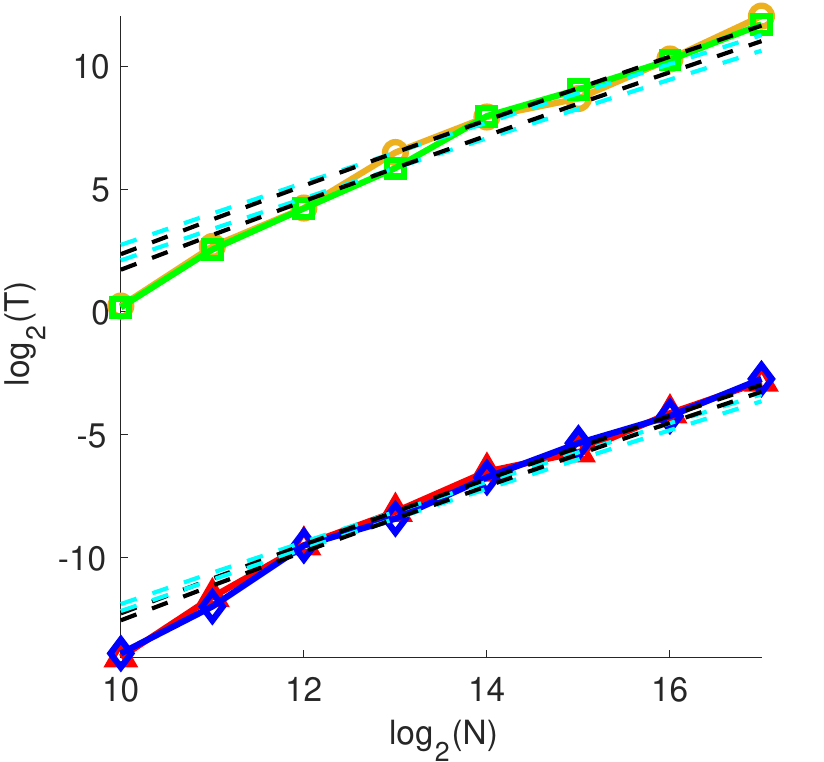} \\
        (3) $m=N$ & (4) $m=1.5N$
    \end{tabular}
    \caption{The computational complexity of the ALT for different orders $m$. Here, $N$ is the size of the matrix, and order $m$ is set to be $0$, $0.5N$, $N$ or $1.5N$. ``Fac'' and ``App'' represent the factorization time and the application time, respectively.  All times are in seconds. {The reference functions $N \Ltwo{N}$ and $N \Lthree{N}$ are multiplied by a constant in order to fit the data.}}
    \label{fig:scaling}
\end{figure}

Figure~\ref{fig:scaling} visualizes the computational complexity of the factorizing and applying the forward and inverse
ALT matrices.  There, $T_{fac}^{fwd}$ and $T_{app}^{fwd}$ are the factorization time and the application time of the 
proposed algorithm for the forward ALT, respectively.  And $T_{mat}^{fwd}$ and $T_{dir}^{fwd}$ are the time for constructing 
the normalized associated Legendre matrix and performing the matrix application directly.
The definitions of $T_{fac}^{inv}$, $T_{app}^{inv}$, $T_{mat}^{inv}$, and $T_{dir}^{inv}$ for the inverse ALT
are analogous.  We observe that the running times of these processes scale nearly linearly with the problem size.

\begin{figure}[!ht]
    \centering
    \begin{tabular}{c}
        \includegraphics[width=11cm,trim={1cm 6.5cm 2.5cm 6.9cm}, clip]{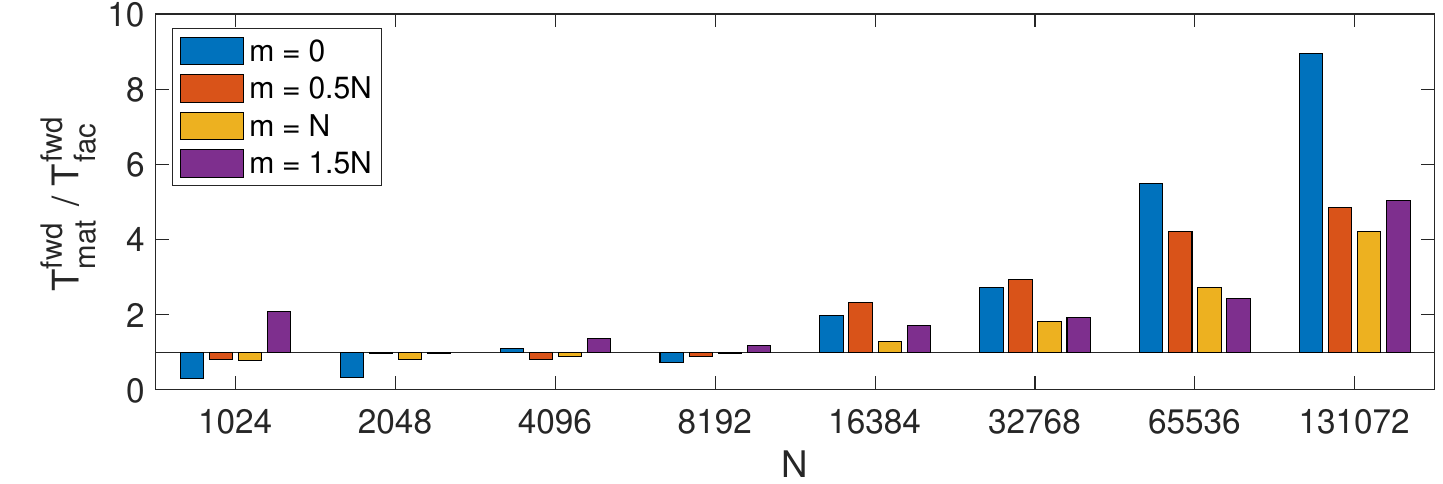} \\
        \includegraphics[width=11cm,trim={1cm 6.5cm 2.5cm 6.9cm}, clip]{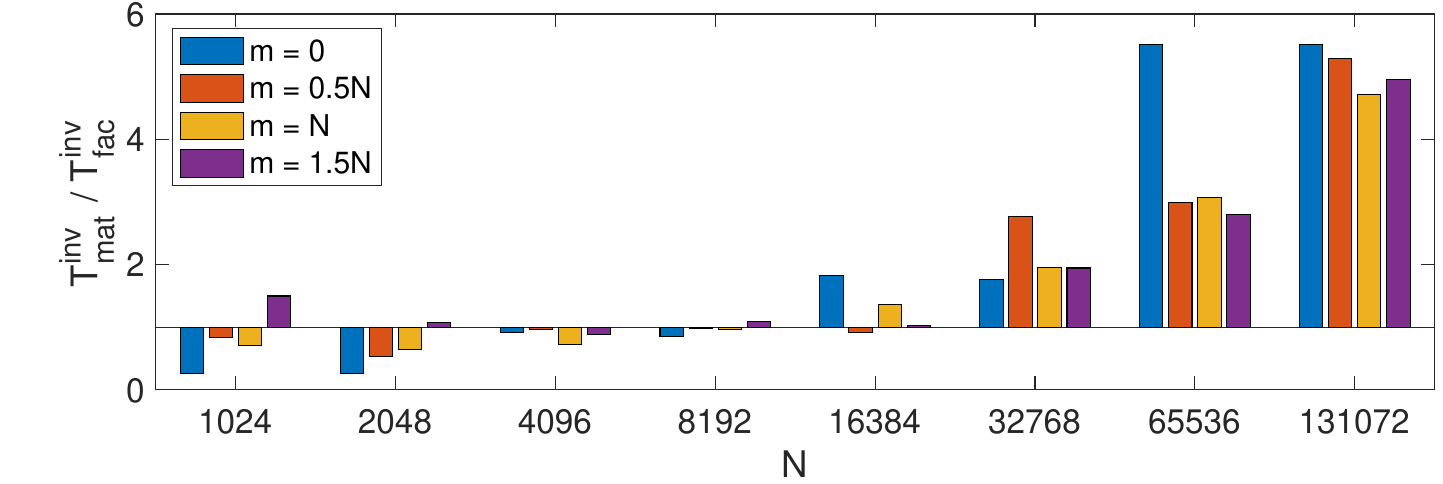} \\
        \includegraphics[width=11cm,trim={1cm 6.5cm 2.5cm 6.9cm}, clip]{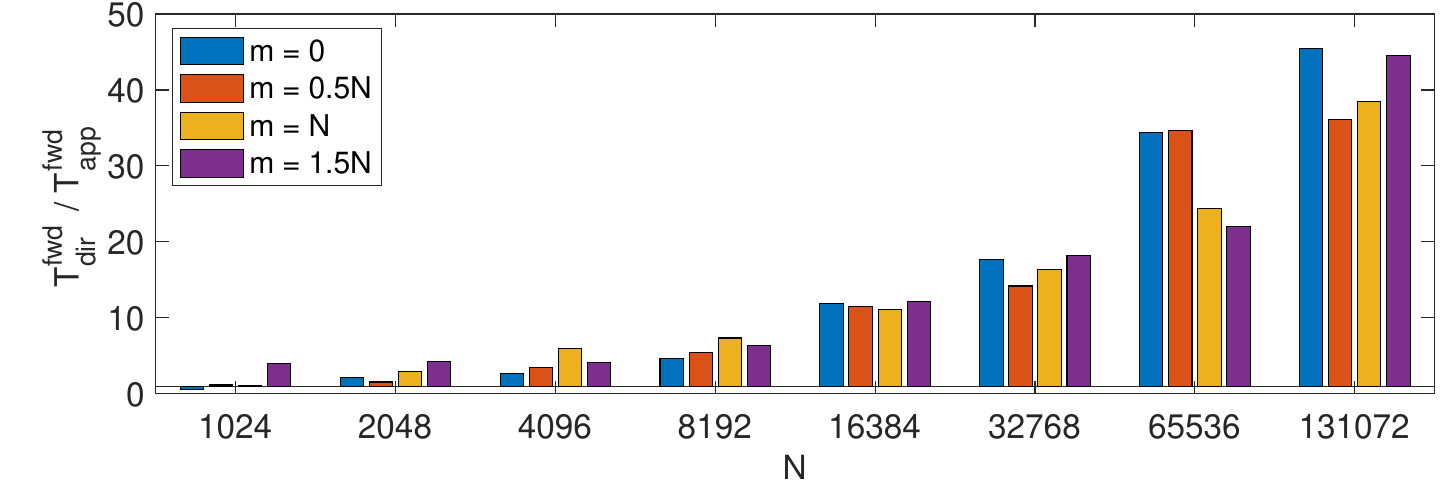} \\
        \includegraphics[width=11cm,trim={1cm 6.5cm 2.5cm 6.9cm}, clip]{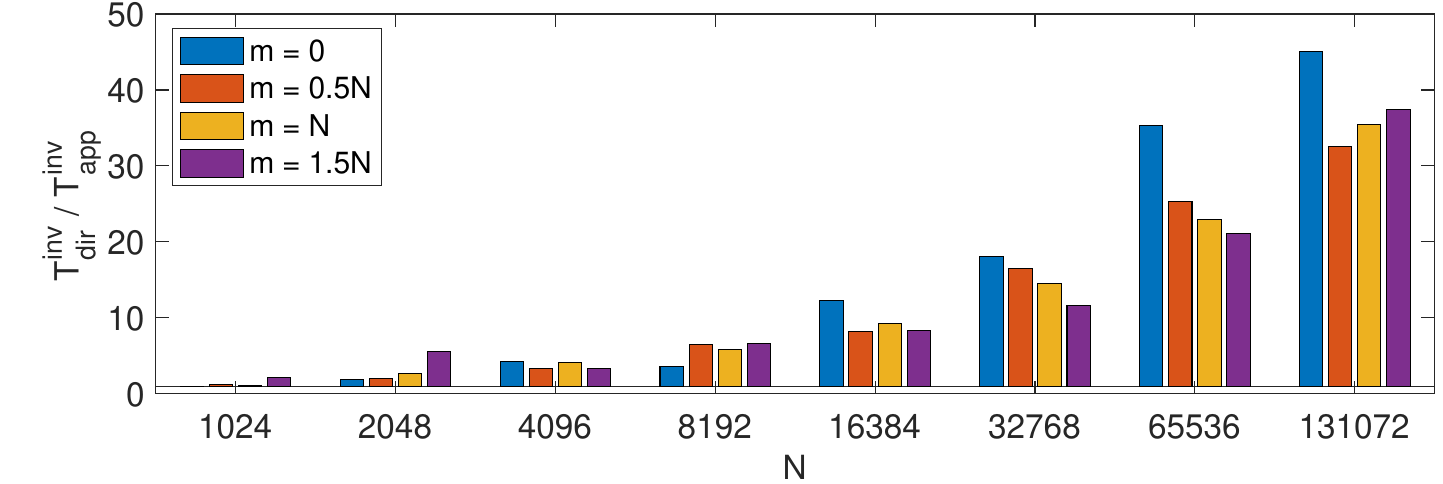}
    \end{tabular}
    \caption{A comparison of the speed of the proposed algorithm for the ALT with the brute force approach.
Here,  $N$ is the size of the matrix, and order $m$ is  $0$, $0.5N$, $N$ or  $1.5N$.
 The adaptive rank $r_k$ for IDBF is set to be $150$.   From the top to bottom, the charts give the
ratios  $T_{mat}^{fwd} / T_{fac}^{fwd}$, $T_{mat}^{inv} / T_{fac}^{inv}$, $T_{dir}^{fwd} / T_{app}^{fwd}$ and 
$T_{dir}^{inv} / T_{app}^{inv}$ are shown.}
    \label{fig:speed}
\end{figure}

Figure~\ref{fig:speed} compares the factorization time and the application time of the proposed algorithm with the brute force approach to applying the ALT (that is, direct application of the matrix discretizing
the ALT).   We observe a significant improvement at larger problem sizes.

\section{Conclusion and future work} \label{sec:conclusion}

This paper introduces an algorithm for the application of the forward and inverse associated Legendre transforms.
Experimental results suggest that its total running time, including both an application and a precomputation phase,
is $\O{N \Lthree{N}}$.  Using this algorithm, the forward and inverse spherical harmonic
transforms can be applied in $\O{N^2 \Lthree{N}}$ time, assuming our conjecture regarding
the running time of our algorithm is correct.

The blocked IDBF algorithm used here is extremely dependent on the method used to form interpolative
decompositions.  The most efficient and accurate method for forming such factorizations is still
an ongoing topic of research, and the authors plan to develop improved versions of their algorithm
which incorporate new developments.

Moreover, the authors are actively working on developing a rigorous bound on the ranks of blocks
of the forward and inverse ALT matrices.  Such a bound would enable a rigorous complexity estimate
for the spherical harmonic transform.

\vspace{0.5cm}
{\bf Acknowledgments.} 
The authors are grateful to the anonymous reviewers for their many helpful comments.
The authors also thank Yingzhou Li for his discussion on block partitioning the oscillatory
region of associated Legendre transform. J.B. was supported in part by NSF grants DMS-1418723 
and DMS-2012487.  Z. C. was partially supported by the Ministry of Education in Singapore under the grant 
MOE2018-T2-2-147. H. Y. was partially supported by NSF under the grant award DMS-1945029.

\bibliographystyle{abbrv}
\bibliography{ref}

\end{document}